\documentclass{aims}
\usepackage{amsmath}
  \usepackage{paralist}
  \usepackage{graphics} 
  \usepackage{epsfig} 
\usepackage{graphicx}  \usepackage{epstopdf}
 \usepackage[colorlinks=true]{hyperref}
\hypersetup{urlcolor=blue, citecolor=red}

\def\currentvolume{X}
 \def\currentissue{X}
  \def\currentyear{200X}
   \def\currentmonth{XX}
    \def\ppages{X--XX}
     \def\DOI{10.3934/xx.xx.xx.xx}

\newtheorem{theorem}{Theorem}[section]
\newtheorem{Theorem}{Theorem}[section]

\newtheorem{Lemma}[theorem]{Lemma}
\newtheorem{Proposition}{Proposition}

\theoremstyle{definition}

\def\disp{\displaystyle}
\def\o{\over}

\def\reals{I\!\!R}
\def\a{\alpha}

\def\g{\gamma}
\def\gm{\gamma(dy,du)}
\def\G{\Gamma}
\def\U{{\mathcal U}}
\def\P{{\mathcal P}}
\def\R{{\mathcal R}}
\def\H{{\mathcal H}}
\def\T{{\mathcal T}}
\def\ph{\varphi}
\def\l{\lambda}
\def\d{\delta}

\def\tt{\Theta}
\def\bu{\bar u}
\def\by{\bar y}

\def\k{\kappa}
\def\b{\beta}
\def\O{\Omega}
\def\bg{\bar \gamma}

\def\hf{\hfill{$\Box$}}

\title[Optimality Conditions and Approximate Solution] 
      {Linear Programming Based Optimality Conditions and Approximate Solution of a Deterministic Infinite Horizon Discounted  Optimal Control Problem in Discrete Time}

\author[Vladimir Gaitsgory, Alex Parkinson and Ilya Shvartsman]{}



\email{vladimir.gaitsgory@mq.edu.au}
 \email{alex.parkinson@students.mq.edu.au}
 \email{ius13@psu.edu}

\thanks{$^*$ Corresponding author}

\begin{document}
\maketitle

\centerline{\scshape Vladimir Gaitsgory}
\medskip
{\footnotesize
 \centerline{Department of Mathematics, Macquarie University, }
   \centerline{ Macquarie Park, NSW 2113, Australia}
} 

\medskip

\centerline{\scshape Alex Parkinson}
\medskip
{\footnotesize
 \centerline{Department of Mathematics, Macquarie University, }
   \centerline{ Macquarie Park, NSW 2113, Australia}
} 

\medskip

\centerline{\scshape Ilya Shvartsman$^*$}
\medskip
{\footnotesize
 \centerline{Department of Mathematics and Computer Science,
Penn State Harrisburg, }
   \centerline{Middletown, PA 17057, USA}
}

\bigskip

 \centerline{(Communicated by the associate editor name)}

\begin{abstract}
It has been recently established that  a deterministic  infinite horizon discounted optimal control problem in discrete time is closely related to a certain infinite dimensional  linear programming  problem and its dual, the latter taking the form of a certain  max-min problem. In the present paper, we use these results to establish necessary and sufficient optimality conditions  for this optimal control problem and to investigate a way  how the latter can be used for the construction of a near optimal control.
\end{abstract}

\section{Introduction and Preliminaries}

It has been  established that a  deterministic infinite horizon discounted optimal control problem in discrete time (for brevity, we will refer to this as just the OC problem) is closely related to an infinite dimensional (ID) linear programming (LP) problem and its dual, the latter taking the form of a certain  max-min problem (see \cite{GPS-1} and also \cite{Jean} for earlier developments in the Markov Decision Processes setting). In the present paper, we use  results of  \cite{GPS-1} to establish necessary and sufficient optimality conditions  for this optimal control problem and to investigate a way  how the latter can be used for the construction of a near optimal control.

Note that necessary optimality conditions in the form of Pontryagin's maximum principle are generally not available for discrete time optimal control problems, this being related to the fact that the so called relaxed (measure valued) controls play no role in dealing with such problems. While the classic control relaxation technique is not applicable in discrete time, a different relaxation approach based on using occupational measures generated by controls and the corresponding solutions of the dynamical system can be used in both continuous and discrete time. Such occupational measures relaxation makes it possible to reformulate the OC problem as an IDLP problem which has a nice dual counterpart, and it is  the relationships between  this IDLP problem and its dual that  lead  to the optimality conditions  established in the paper.

The linear programming (LP) approach to optimal control problems has been studied extensively in both stochastic and deterministic settings (see, e.g.,  \cite{BhBo}, \cite{Vivek}, \cite{BGQ}, \cite{F-V},\cite{Her-Her-Lasserre},  \cite{Kurtz},  \cite{Stockbridge}, \cite{Stockbridge1} and, respectively, \cite{Adelman}, \cite{Gai-F-Leb}, \cite{Gai8}, \cite{GQ}, \cite{GQ-1}, \cite{GR}, \cite{Goreac-Serea},  \cite{Kamo}, \cite{Adelman-1},  \cite{Lass-Trelat}, \cite{QS}, \cite{Rubio}, \cite{Vinter} as well as references therein). In particular, results establishing the validity of LP formulations of  deterministic infinite horizon OC problems with time discounting  have been obtained in \cite{GQ}, \cite{GQ-1} and \cite{Kamo} for systems evolving in continuous time and in \cite{GPS-0} and \cite{GPS-1}  for systems evolving in discrete time. (Note that other approaches/techniques for dealing with deterministic optimal control problems on the infinite time horizon have been studied, e.g., in   \cite{BCD}, \cite{B12}, \cite{CHL}, \cite{GruneSIAM98}, \cite{GruneJDE98}, \cite{Z14}, \cite{Z06} and \cite{Z06a}; see also references therein.) In this paper, we continue the line of research started in  \cite{GPS-0} and \cite{GPS-1}.

 Consider the control  system
\begin{equation}\label{A1}
\begin{aligned}
&y(t+1)=f(y(t),u(t)),\ \; t\in \T :=\{0,1,\dots\ \}, \\
&y(0)=y_0,\\
&y(t)\in Y,\\
&u(t)\in U(y(t)),
\end{aligned}
\end{equation}
where  $Y$ is a given nonempty compact subset of $\reals^m$, $\ U(\cdot):\,Y\leadsto U_0$ is an upper semicontinuous compact-valued mapping to a given compact metric space $U_0$,
$\ f(\cdot, \cdot ):\,\reals^m\times U_0\to \reals^m$ is a continuous function.
  A control $u(\cdot)$ and the pair $(y(\cdot),u(\cdot))$ are called an admissible control and, respectively, an admissible process if the relationships \eqref{A1} are satisfied. The set of all admissible controls is denoted as $\U(y_0)$.

 Everywhere in what follows, we will be
 dealing with the optimal control problem
\begin{equation}\label{A111}
\inf_{u(\cdot)\in \U(y_0)}\sum_{t=0}^{\infty} \a^t g(y(t),u(t))=:V(y_0),
\end{equation}
where
$g:\,\reals^m\times U_0\to \reals^m$ is a continuous function and  $\a\in (0,1)$ is a discount factor.
Note that the last two constraints in \eqref{A1} can be rewritten as one:
$$
u(t)\in A(y(t)),
$$
where the map $\ A(\cdot):\,Y\leadsto U_0$ is defined by the equation
\begin{equation}\label{e-A}
\begin{aligned}
A(y):=\{u\in U(y)|\, f(y,u)\in Y\} \ \ \ \forall y\in Y.
\end{aligned}
\end{equation}
Note also that from the fact that $ U(\cdot)$ is upper semicontinuous and $\ f(\cdot, \cdot )$ is continuous it follows that the map $A(\cdot)$ is upper semicontinuous and its graph $G$,
$$
G:={\rm graph}\,A=\{(y,u)|\,y\in Y,\,u\in U(y),\,f(y,u)\in Y\},
$$
is a compact subset of $Y\times U_0$.

 The standing assumption in the paper is the following

  ASSUMPTION I: {\em The set $A(y)$ is not empty for any $y\in Y$.}

 As can be readily seen,  this assumption implies that the set  $\U(y)$ is not empty for any $y\in Y$  (systems  satisfying such a property are called {\em viable}; see \cite{Aub}). It can be shown  (see, e.g., \cite{Bert} or \cite{GPS-1}) that under this assumption an optimal solution of problem \eqref{A1} exists, the optimal value function $V(\cdot)$ is lower semicontinuous and bounded on $Y$, and $V(\cdot)$ is a solution of the  following equation   (the dynamic programming principle)
  \begin{equation}\label{B1}
V(y)=\min_{u\in A(y)}\{g(y,u)+\a V(f(y,u))\}\ \ \ \ \forall y\in Y.
\end{equation}
Note that this equation can be rewritten in the  form that resembles the Hamilton-Jacobi-Bellman equation for continuous time systems:
\begin{equation}\label{E3}
H_{V}(y)-(1-\a)V(y)=0\ \ \ \ \ \forall y\in Y,
\end{equation}
where, for  any lower semicontinuous function $\psi:\,Y\to \reals$,  
\begin{equation}\label{e-HJB}
\begin{aligned}
H_{\psi}(y):=\min_{u\in A(y)}\{g(y,u) + \a (\psi(f(y,u))-\psi(y))\} .
\end{aligned}
\end{equation}
Along with (\ref{E3}), let us consider the following max-min problem
\begin{equation}\label{B5}
\max_{\psi\in LS}\inf_{y\in Y}\{H_{\psi}(y)-(1-\a)(\psi(y)-\psi(y_0))\}=:\mu^*(y_0)
\end{equation}
%
where $max$ is over the class of bounded lower semicontinuous functions from $Y$ to $\reals$ (denoted as $LS$). It has been recently established in \cite{GPS-1} (see also Proposition \ref{P5} below) that the maximum in \eqref{B5} is reached at $\psi=V$.
Note that from the fact that $V$ is a maximizer in \eqref{B5} it follows that $\tilde V:=V+{\rm const}$ is a maximizer in \eqref{B5} as well.
The set of maximizers in \eqref{B5} may, in fact, be much broader than just constant shifts of the optimal value function (see an example in Section \ref{Section-NSOC}).

               In this paper, we  establish that necessary and sufficient optimality conditions for problem (\ref{A111})  can be stated in terms of any such maximizer. (Note that the max-min problem (\ref{B5}) involves the dependence on $y_0$ and the optimality conditions in terms of a maximizer in (\ref{B5}) will only be  valid for the solutions satisfying the initial condition $y(0)=y_0 $, this being  in contrast to the solution of the dynamic programming equation (\ref{E3}), which allows to characterize the optimal solutions for arbitrary initial conditions.)
               
               We  will also indicate a way how an approximate maximizer in  \eqref{B5}  can be used for
              the construction of a near optimal control for problem
(\ref{A111}), this construction will be illustrated by a numerical example.

The paper is organized as follows.
In Section \ref{Section-NSOC}, we establish necessary and sufficient conditions of optimality for problem (\ref{A111}) in terms of a maximizer in the max-min problem (\ref{B5}) (the main result of this section is Theorem \ref{P2.1}). In Section \ref{Section-approx-max-min}, we introduce  $N$-approximating max-min problem, in which
 in contrast to  (\ref{B5}),  the maximization  is over $\psi\in D_N$, where $D_N$ is an $N$-dimensional subspace of the space of continuous functions on $Y$. We show that the optimal value of the $N$-approximating max-min problem converges to the optimal value of
 (\ref{B5}) as $N\rightarrow\infty $ (Proposition \ref{P6}) and we establish that a maximizer in the $N$-approximating max-min problem exists for any $N$ under a readily verifiable controllability condition (Proposition \ref{P3.2}). In Section \ref{Section-construction},
 we establish that a maximizer in the $N$-approximating max-min problem can be used for the construction of a near optimal control (Theorem \ref{T1}). Sections \ref{Section-example} and \ref{heuristic} are devoted to numerical examples that illustrate the latter construction. Some of the results presented in
 Section \ref{Section-NSOC} were  announced in \cite{GPS-0} without proofs. Continuous time counterparts of  results of Sections \ref{Section-approx-max-min} and \ref{Section-construction} can be found in \cite{GRT}.

To conclude this section, let us outline some notations and results that are used further in the text.
For  an admissible process $(y(\cdot),u(\cdot))$, a probability measure $\g_{u}$ is called the {
\em discounted occupational measure} generated by $u(\cdot)$ if, for any Borel set $Q\subset G$,
\begin{equation}\label{EE6}
\g_{u}(Q)=(1-\a) \sum_{t=0}^{\infty} \a^t 1_Q(y(t),u(t)),
\end{equation}
where $1_Q(\cdot)$ is the indicator function of $Q$. 

It can be shown that this definition is equivalent to the validity of the relationship
\begin{equation}\label{G8}
\int_{G} q(y,u) \g_u(dy,du)=(1-\a)\sum_{t=0}^{\infty} \a^t q(y(t),u(t))
\end{equation}
for any Borel measurable function $q$ on $G$. Indeed, \eqref{EE6} obviously implies the validity of \label{G8} for a simple function (i.e., a finite sum of indicator functions of Borel measurable sets).  The validity of \label{G8} for an arbitrary Borel $q$  follows from the definition of the Lebesgue integral as a limit of integrals of simple functions; see, e.g. \cite{Ash}. 

To describe convergence properties of occupational measures, we introduce the following metric on $\P(G)$ (the space of probability measures defined on Borel subsets of $G$):
$$
\rho(\g',\g''):=\sum_{j=1}^{\infty} {1\o 2^j}\left|\int_G q_j(y,u)\g'(dy,du)-\int_G q_j(y,u)\g''(dy,du)\right|
$$
for $\g',\g''\in \P(G)$, where $q_j(\cdot),\,j=1,2,\dots,$ is a sequence of Lipschitz continuous functions dense in the unit ball of the space of continuous functions $C(G)$ from $G$ to $\reals$.
This metric is consistent with the weak$^*$ convergence topology on $\P(G)$, that is,
a sequence $\g^k\in \P(G)$ converges to $\g\in \P(G)$ in this metric if and only if
$$
\lim_{k\to \infty}\int_G q(y,u)\g^k(dy,du)=\int_G q(y,u)\g(dy,du)
$$
for any $q\in C(G)$. Using this metric, we can define the ``distance" $\rho(\g,\Gamma)$ between $\g\in \P(G)$ and $\Gamma\subset \P(G)$
and the Hausdorff metric $\rho_H(\Gamma_1,\Gamma_2)$ between $\Gamma_1\subset \P(G)$ and $\Gamma_2\subset \P(G)$ as follows:
$$
\rho(\g,\Gamma):=\inf_{\g'\in \Gamma}\rho(\g,\g'),\quad
\rho_H(\Gamma_1,\Gamma_2):=\max\{\sup_{\g\in \Gamma_1}\rho(\g,\Gamma),\sup_{\g\in \Gamma_2}\rho(\g,\Gamma_2)\}.
$$
Let $\G(y_0)$ denote the set of  all discounted occupational measures generated by the admissible controls. That is,
$$
\G(y_0):=\bigcup_{u(\cdot)\in \U}\{\g_{u}\}.
$$
Note that, due to Assumption I, the set $\G(y_0) $ is not empty. Also, due to \eqref{G8},  problem \eqref{A1}
can be rewritten as
\begin{equation}\label{A3}
\min_{\g\in \G(y_0)} \int_{G} g(y,u)\g(dy,du)= (1-\a)V(y_0).
\end{equation}


Consider the problem
\begin{equation}\label{D1}
\min_{\g\in W(y_0)} \int_{G} g(y,u)\g(dy,du)=:g^*(y_0),
\end{equation}
where
$W(y_0)$  is  the set of probability measures defined by the equation
\begin{equation}\label{e-W-y-0}
\begin{aligned}
&W(y_0):=\{\g\in \P(G)|\, \\
&\int_{G}[\a(\ph(f(y,u))-\ph(y))
+(1-\a)(\ph(y_0)-\ph(y))]\gm=0\quad \forall \ph\in C(Y)\}.
\end{aligned}
\end{equation}
Note that (\ref{D1}) is an infinite dimensional linear programming  problem since both the objective function and the constraints are linear with respect to the \lq\lq decision variable" $\gamma$. Note also that the set $W(y_0) $ is not empty (as follows from the statement (i) of Proposition \ref{P5} below) and, as can be readily verified, it is compact in weak$^*$ topology. Hence, the minimum in (\ref{D1}) is reached.

In \cite{GPS-1}, it has been established that problem (\ref{B5}) is dual  to the IDLP problem (\ref{D1}) and that these two problems  are related to the optimal control problem (\ref{A111}). Some of the  relationships between problems (\ref{B5}), (\ref{D1}) and (\ref{A111})  are summarized in the following proposition.

\begin{Proposition}\label{P5} The following statements are valid:

(i) The closed convex hull of the set of discounted occupational measures $\G(y_0)$ is equal to the set $W(y_0)$. That is,
\begin{equation}\label{e-disc-equality}
\bar{\rm co}\,\G(y_0)=W(y_0).
\end{equation}

(ii) The optimal values in problems \eqref{B5} and \eqref{D1} coincide and are equal to the optimal value in (\ref{A111}) multiplied by $(1-\a)$, that is,
\begin{equation}\label{D3}
\mu^*(y_0)=g^*(y_0)=(1-\a)V(y_0).
\end{equation}

(iii) The supremum in \eqref{B5} is reached at $\psi=V$.

\end{Proposition}

 {\bf Proof.} The statements (ii) and (iii) follow from Theorem 4.1 in \cite{GPS-1}. The statement (i) follows from Corollary 2 in \cite{GPS-1}.
\hf

\section{Necessary and Sufficient Optimality Conditions}\label{Section-NSOC}

We will say that $ \psi\in LS$ is a solution of \eqref{B5} if
\begin{equation}\label{B5-2}
\inf_{y\in Y}\{H_{\psi}(y)-(1-\a)(\psi(y)-\psi(y_0))\}=\mu^*(y_0),
\end{equation}
which is equivalent to
\begin{equation}\label{e-opt-dual-solution}
\inf_{(y,u)\in G}\{g(y,u)+\a(\psi(f(y,u))-\psi(y))+(1-\a)(\psi(y_0)-\psi(y))\}=\mu^*(y_0).
\end{equation}
Due to Proposition \ref{P5}(iii), $V$ is a maximizer in \eqref{B5}, and, as has been mentioned earlier, all constant shifts of $V$ are solutions of \eqref{B5} as well. As also has been mentioned above, the set of solutions of problem \eqref{B5} can be significantly larger than the set of these shifts. This is demonstrated by the following example.

{\bf Example.} Consider the problem
\begin{equation*}
\begin{aligned}
&{\rm Minimize}\,\sum_{t=0}^{\infty} \a^t g(y(t)),\\
&y(t+1)=u(t),\,t\in \T,\\
&y(0)=y_0,\\
&u(t)\in [0,1],\\
&y(t)\in [0,1],
\end{aligned}
\end{equation*}
where function $g$ is increasing on $[0,1]$ and $g(0)=0$.

It is clear that the optimal control is $u\equiv 0$ with the corresponding trajectory
\begin{equation*}
y(t)=\begin{cases}
y_0,&t=0,\\
0, &t>0,
\end{cases}
\end{equation*}
and the value function is $V(y)=g(y)$.

Let us show that, if $\psi:\,[0,1]\to \reals$ is such that $\psi(y_0)=g(y_0)$, $\psi(0)=g(0)=0$, and $0\le\psi(y)\le g(y)$ for all $y\in [0,1]$, then
$\psi$ is a solution of \eqref{B5}. Indeed, for such $\psi$ we have
$$
\min_{(y,u)\in [0,1]\times[0,1]}\{g(y)+\a\psi(u)-\psi(y)+(1-\a)\psi(y_0)\}=(1-\a)\psi(y_0)=(1-\a)V(y_0),
$$
therefore, $\psi$ is a solution of \eqref{B5} due to Proposition \ref{P5}(ii). \hf

\bigskip


The theorem stated below establishes necessary and sufficient optimality conditions for problem (\ref{A111}) in terms of any solution of \eqref{B5}.

\begin{Theorem}\label{P2.1}
Let $\psi$ be a solution of \eqref{B5}.
Optimality of an admissible process $(y(\cdot),u(\cdot))$ is equivalent to each of the following for all $t\in \T$:

{\rm (i)}
\begin{equation}\label{N3}
(y(t),u(t))={\rm argmin}_{(y,u)\in G}\{g(y,u)+\a\psi(f(y,u))-\psi(y)\},
\end{equation}
or, equivalently,
\begin{equation}\label{A16}
\begin{aligned}
&u(t)={\rm argmin}_{u\in A(y)}\{g(y(t),u)+\a\psi(f(y(t),u))\},\\
&y(t)={\rm argmin}_{y\in Y}\{H_{\psi}(y)-(1-\a)\psi(y)\};
\end{aligned}
\end{equation}
{\rm (ii)}
\begin{equation}\label{A17}
\psi(y(t))=V(y(t))+\psi(y_0)-V(y_0),
\end{equation}
that is, $\psi$ and $V$ coincide on the optimal trajectory up to the constant $\psi(y_0)-V(y_0)$.

Furthermore, if  $(y(\cdot),u(\cdot))$ is optimal, then
\begin{equation}\label{A9}
H_{\psi}(y(t))-(1-\a)\psi(y(t))=(1-\a)(V(y_0)-\psi(y_0)).
\end{equation}
\end{Theorem}

{\bf Proof.} Since both $V$ and $\psi$ are solutions of \eqref{B5}, we have
\begin{equation*}
\inf_{y\in Y}\{H_{V}(y)+(1-\a)(V(y_0)-V(y))\}=
\inf_{y\in Y}\{H_{\psi}(y)+(1-\a)(\psi(y_0)-\psi(y))\}.
\end{equation*}
Taking into account that $H_{V}(y)-(1-\a)V(y)= 0$ for all $y\in Y$ due to \eqref{E3}, we have
\begin{equation}\label{N4}
\inf_{y\in Y}\{H_{\psi}(y)-(1-\a)\psi(y)\}=(1-\a)(V(y_0)-\psi(y_0)),
\end{equation}
which implies that
\begin{equation}\label{C2}
g(y,u)+\a\psi(f(y,u))-\psi(y)\ge (1-\a)(V(y_0)-\psi(y_0)) \quad \hbox{for all }(y,u)\in G.
\end{equation}
Further, the identity
\begin{equation}\label{N1}
\sum_{t=0}^{\infty} \a^t \psi(y(t))=\psi(y_0)+\a\sum_{t=0}^{\infty} \a^t \psi(f(y(t),u(t))),
\end{equation}
is equivalent to
\begin{equation}\label{E4}
\sum_{t=0}^{\infty} \a^t [\a\psi(f(y(t),u(t)))-\psi(y(t))+(1-\a)\psi(y_0)]=0.
\end{equation}
Assume that $(y(\cdot),u(\cdot))$ is an optimal process. Then
$$
\sum_{t=0}^{\infty} \a^t g(y(t),u(t))=V(y_0),
$$
or, equivalently,
\begin{equation}\label{A15}
\sum_{t=0}^{\infty} \a^t [g(y(t),u(t))-(1-\a)V(y_0)]=0.
\end{equation}
Adding \eqref{E4} and \eqref{A15} we obtain
\begin{equation}\label{A12}
\sum_{t=0}^{\infty} \a^t [g(y(t),u(t))+\a\psi(f(y(t),u(t))-\psi(y(t))+(1-\a)(\psi(y_0)-V(y_0))]=0.
\end{equation}
Taking into account \eqref{C2} we arrive at
\begin{equation}\label{C4}
g(y(t),u(t))+\a\psi(f(y(t),u(t)))-\psi(y(t))= (1-\a)(V(y_0)-\psi(y_0)) \quad \hbox{for all }t,
\end{equation}
and \eqref{N3}.

Let us see that  \eqref{A9} holds.
From the definition of $H$, \eqref{A9} is equivalent to
\begin{equation*}
\min_{u\in A(y(t))}\{g(y(t)),u)+\a\psi(f(y(t),u))\}-\psi(y(t))= (1-\a)(V(y_0)-\psi(y_0)) \quad \hbox{for all }t.
\end{equation*}
From \eqref{C2} and \eqref{C4} it follows that this relation holds with the minimizing $u\in A(y(t))$ equal to $u(t)$.

To prove \eqref{A17} take $k\in \T$. From the identity
\begin{equation*}
\sum_{t=k}^{\infty} \a^t \psi(y(t))=\a^k\psi(y(k))+\sum_{t=k+1}^{\infty} \a^t \psi(y(t))
\end{equation*}
we get
\begin{equation}\label{H1}
\a^k\psi(y(k))=\sum_{t=k}^{\infty} \a^t [\psi(y(t))-\a\psi(f(y(t),u(t))].
\end{equation}
From optimality of $(y(\cdot),u(\cdot))$, via \eqref{C4}, we obtain
\begin{equation}\label{H2}
\begin{aligned}
\a^k\psi(y(k))&=\sum_{t=k}^{\infty} \a^t[g(y(t),u(t))+(1-\a)(\psi(y_0)-V(y_0))]\\
&=\a^k\left\{\sum_{t=0}^{\infty} \a^t g(y(t+k),u(t+k))+(\psi(y_0)-V(y_0))\right\}\\
&=\a^k \{V(y(k))+\psi(y_0)-V(y_0)\},
\end{aligned}
\end{equation}
therefore, $\psi(y(k))=V(y(k))+V(y_0)-\psi(y_0)$, i.e., \eqref{A17} holds.

Conversely, assume that assertion (i) of the theorem is true. Let us rewrite \eqref{N4} in the form
$$
\inf_{y\in Y}\min_{u\in A(y)}\{g(y,u)+\a\psi(f(y,u))-\psi(y)\}=(1-\a)(V(y_0)-\psi(y_0)).
$$
From \eqref{N3} it follows that $\disp\inf_{y\in Y}\min_{u\in A(y)}$ in this formula is reached when $(y,u)=(y(t),u(t))$ for all $t$,  in which case \eqref{C4} holds. The latter
implies \eqref{A12}. Subtracting \eqref{E4} from \eqref{A12} we obtain \eqref{A15}, therefore, the process $(y(\cdot),u(\cdot))$ is optimal.

Assume now that (ii) holds. Similarly to formula \eqref{H2}, we obtain for any $k\in \T$
\begin{equation*}
\begin{aligned}
&\psi(y(k))=V(y(k))+\psi(y_0)-V(y_0)=\sum_{t=0}^{\infty} \a^t g(y(t+k),u(t+k))+(\psi(y_0)-V(y_0))\\
&=\a^{-k}\left\{\sum_{t=k}^{\infty} \a^t[g(y(t),u(t))+(1-\a)(\psi(y_0)-V(y_0))]\right\}.
\end{aligned}
\end{equation*}
From this equality and \eqref{H1} we obtain
\begin{equation*}
\sum_{t=k}^{\infty} \a^t [g(y(t),u(t))+(1-\a)(\psi(y_0)-V(y_0))]=\sum_{t=k}^{\infty} \a^t[\psi(y(t))-\a\psi(f(y(t),u(t))],
\end{equation*}
wherein setting $k=0$ leads to \eqref{A12}.
As above, subtracting \eqref{E4} from \eqref{A12} we get \eqref{A15}, that is,  the process $(y(\cdot),u(\cdot))$ is optimal. The theorem is proved.
\hf

\bigskip

Assertion (ii) of the theorem above states that $\psi$ and $V$ are equal on the
optimal trajectory, up to the constant $\psi(y_0)-V(y_0)$. It can be shown
that, away from the optimal trajectory, the equality becomes inequality.
\begin{Proposition}\label{P7}
Let $\psi$ be a solution of \eqref{B5}. Then $\psi(y)\le V(y)+\psi(y_0)-V(y_0)
$.
\end{Proposition}
{\bf Proof.} Take the optimal trajectory $(y(\cdot),u(\cdot))$ starting
from an arbitrary initial condition $y_1$.
Replacing $y_0$ with $y_1$ in \eqref{N1} we get
\begin{equation}\label{N2}
\psi(y_1)=\sum_{t=0}^{\infty} \a^t(\psi(y(t))-\a\psi(f(y(t),u(t)))).
\end{equation}
From \eqref{C2} it follows that
$$
\psi(y(t))-\a\psi(f(y(t),u(t))\le g(y(t),u(t))-(1-\a)(V(y_0)-\psi(y_0)).
$$
Taking into account \eqref{N2} and the latter inequality we obtain
$$
\psi(y_1)\le \sum_{t=0}^{\infty} \a^t g(y(t),u(t))-(V(y_0)-\psi(y_0))= 
V(y_1)-(V(y_0)-\psi(y_0)),
$$
and the assertion of the proposition follows.
\hf

\bigskip

{\bf Remark.} If we introduce the function resembling the unmaximized Hamiltonian
$$
\H_{\psi}(y,u):=\a(\psi(f(y,u))-\psi(y))+g(y,u),
$$
then the first of the conditions \eqref{A16} can be written in the form resembling the Pontryagin-type minimum principle
$$
u(t)={\rm argmin}_{u\in A(y)}\,\H_{\psi}(y(t),u).
$$
(Of course, it is well known that the Pontryagin maximum principle doesn't hold in general for discrete time systems without additional convexity assumptions.)

\bigskip

Let us complete this section with  another characterization of the solution set of the max-min problem (\ref{B5}). Consider the inequality (compare with (\ref{E3}))
\begin{equation}\label{E3-inequality-1}
H_{\psi}(y)-(1-\a)\psi(y)\geq 0 \ \ \ \ \ \forall y\in Y.
\end{equation}
The following result establishes the relationships between the solutions  of this inequality
satisfying the additional  condition
\begin{equation}\label{E3-in-condition}
\psi(y_0)= V(y_0)
\end{equation}
and the solutions of the max-min problem (\ref{B5}).

\begin{Proposition}\label{P5-1}
If $\psi\in LS$ is a solution of inequality (\ref{E3-inequality-1}) satisfying
(\ref{E3-in-condition}), then $\psi$ is a solution of the max-min problem \eqref{B5}. Conversely, if  $\psi$ is a solution of the max-min problem \eqref{B5}, then 
\begin{equation}\label{e-disc-equivalence}
\tilde \psi (y):= \psi(y) - \psi(y_0) + V(y_0),
\end{equation}
is a solution of inequality (\ref{E3-inequality-1}) satisfying
(\ref{E3-in-condition}).
\end{Proposition}

{\bf Proof.} Let $\psi\in LS$ be a solution of inequality (\ref{E3-inequality-1}) satisfying
(\ref{E3-in-condition}). Then, by (\ref{D3}),
$$
H_{\psi}(y)-(1-\a)(\psi(y)-\psi(y_0))\geq (1-\a)\psi(y_0)= (1-\a)V(y_0)=\mu^*(y_0)  \ \ \ \ \ \forall y\in Y.
$$
The latter implies (\ref{B5-2}) due to the definition of $\mu^*(y_0)$ (see (\ref{B5})). Hence, $\psi$ is a solution of the max-min problem \eqref{B5}.

Let now $\psi$ be a solution of the max-min problem \eqref{B5}. Then $\tilde \psi $ defined in (\ref{e-disc-equivalence})  will be a solution of this problem too. Also, $\tilde \psi $ will satisfy  (\ref{E3-in-condition}) ($\tilde \psi(y_0) = V(y_0) $). Consequently,
by (\ref{D3}),
$$
H_{\tilde \psi}(y)-(1-\a)(\tilde \psi(y)-\tilde \psi(y_0))\geq \mu^*(y_0)= (1-\a)V(y_0) = (1-\a)\tilde \psi(y_0)
$$
$$
\Rightarrow \ \ \ \ \ \ H_{\tilde \psi}(y)-(1-\a)\tilde \psi(y))\geq 0 \ \ \ \ \ \forall y\in Y.
$$
That is, $\tilde \psi $ is a solution of inequality (\ref{E3-inequality-1}) as well.
\hf

\section{$N$-Approximating Max-Min Problem}\label{Section-approx-max-min}


Let $\{\phi_i\}_{i=1}^{\infty}$ be a sequence of continuous functions with the following properties:

(i) Any finite collection of functions from this sequence is linearly independent on any open set;

(ii) For any $\psi\in C(Y)$ (the space of continuous functions on $Y$) and any $\d>0$, there exist $N$ and $\l_i^N$, $i=1,\dots,N$ such that $\disp \max_{y\in Y}|\psi(y)-\sum_{i=1}^N \l_i^N\phi_i(y)|\le \d$.
\medskip

Note that, due to the \lq\lq approximation property" (ii), the set $W(y_0)$ defined in (\ref{e-W-y-0}) can be rewritten in the form
\begin{equation}\label{e-W-y-0-countable}
\begin{aligned}
&W(y_0):=\{\g\in \P(G)|\, \\
&\int_{G}[\a(\ph_i(f(y,u))-\ph_i(y))
+(1-\a)(\ph_i(y_0)-\ph_i(y))]\gm=0\quad \forall \ i=1,...\}.
\end{aligned}
\end{equation}

\medskip

 In what follows it is assumed that $\phi_1\equiv 1$. Hence, the linear independence property (i) implies that
\begin{equation}\label{E2}
\sum_{i=1}^N v_i \phi_i(y) =\hbox{const}\ \   \ \forall y\in Q\ \ \ \ \ \Leftrightarrow\ \ \ \
v_i =0 \ \forall\ i=1,\dots,N
\end{equation}
if $Q$ has a nonempty interior. (An example of a sequence with the properties (i) and (ii) is the sequence of monomials $y_1^{i_1}\dots y_m^{i_m},\,i_1,\dots,i_m=0,1,\dots$, where $y_j$ stands for the $j$th component of $y$.)

Consider the max-min problem
\begin{equation}\label{eq-maxmin-N}
\begin{aligned}
\sup_{\psi\in D_N}\min_{(y,u)\in G}\{g(y,u)+\a(\psi(f(y,u))-\psi(y))+(1-\a)(\psi(y_0)-\psi(y))\}=:\mu_N^*(y_0),
\end{aligned}
\end{equation}
where $D_N\subset C(Y)$ is the finite dimensional space defined by the equation
\begin{equation*}
D_N:=\{\psi\in C(Y)|\, \psi(y)=\sum_{i=1}^N \l_i \phi_i(y);\ \l_i\in \reals,\, i=1,\dots,N\}.
\end{equation*}
Problem (\ref{eq-maxmin-N})
is referred to as the {\em $N$-approximating problem}.


\bigskip
\begin{Proposition}\label{P6} The optimal value  $\mu_N^*(y_0) $ of the $N$-approximating problem converges to the optimal value $\mu^*(y_0)$ of problem \eqref{B5} as $N$ tends to infinity. That is,
\begin{equation*}\label{D6}
\lim_{N\to \infty}\mu^*_N(y_0)=\mu^*(y_0).
\end{equation*}
\end{Proposition}

{\bf Proof.}
Let
\begin{equation}\label{D10}
\mu_C^*(y_0):=\sup_{\psi\in C(Y)}\min_{(y,u)\in G}\{g(y,u)+\a(\psi(f(y,u))-\psi(y))+(1-\a)(\psi(y_0)-\psi(y))\},
\end{equation}
where
 {\it sup} is over the space of continuous functions.
 It is obvious that
  $\mu^*_N(y_0)\leq \mu_C^*(y_0)\ \forall N$. Also  
the sequence $\{\mu^*_N(y_0)\}$ is monotone increasing. Hence, the limit $\disp \lim_{N\to \infty}\mu^*_N(y_0)$ exists and does not exceed $\mu_C^*(y_0)$. Also, from the approximating property of $\{\phi_i\}$ it follows, in fact,  that $\disp \lim_{N\to \infty}\mu^*_N(y_0)=\mu_C^*(y_0)$.
 Thus, to prove the proposition, it is sufficient to establish that $\mu_C^*(y_0) $ is equal to $\mu^*(y_0)$. This is established by the lemma stated below. \hfill{$\Box$}

 \begin{Lemma}\label{e-duality-inf}
 The equality
\begin{equation}\label{e:muC=mu}
\mu_C^*(y_0)=\mu^*(y_0)
\end{equation}
is valid.
 \end{Lemma}

{\bf Proof.}  Since
\begin{equation}\label{e-inequality}
\mu_C^*(y_0)\leq\mu^*(y_0),
\end{equation}
it is sufficient to establish that the inequality  $\mu_C^*(y_0)\geq\mu^*(y_0) $ holds true.

To prove the latter,  note first that for any $\psi\in C(Y)$, we have
\begin{equation}\label{H3}
\min_{(y,u)\in G}\{\a(\psi(f(y,u))-\psi(y))+(1-\a)(\psi(y_0)-\psi(y))\}\le 0.
\end{equation}
Indeed, if this was not the case, then, for $\psi_m:=m\psi$ with positive integer $m$ we would get
$$
\lim_{m\to \infty}\min_{(y,u)\in G}\{g(y,u)+\a(\psi_m(f(y,u))-\psi_m(y))+(1-\a)(\psi_m(y_0)-\psi_m(y))\}=+\infty,
$$
contradicting the fact that $\mu_C^*$  is bounded (the latter being implied by (\ref{e-inequality})).

Assume that functions $\{\phi_i\}$ are normalized so that $\max_{y\in Y}|\phi_i(y)|<1/2^i$. Define $\hat{Q}\subset\reals\times l_1$ by
		\begin{equation*}
				\begin{aligned}
			&\hat{Q} := \{(\theta,x) |\, \theta\geq\int_{G}g(y,u)\gamma(dy,du),\, x=\left(x_1,x_2,\dots\right),  \\
			&x_i=\int_{G}[\alpha(\phi_i\left(f\left(y,u\right)\right)-\phi_i\left(y\right)) + \left(1-\alpha\right)(\phi_i\left(y_0\right)-\phi_i\left(y\right))]\gamma\left(dy,du\right),\,\gamma\in\mathcal{P}(G)\}.
		\end{aligned}
		\end{equation*}
It's easy to see that the set $\hat{Q}$ is a closed subset of $\reals\times l_1$ and that, for any $j=1,2,\dots$, the point $(g^*(y_0)-\frac{1}{j},0)$ does not belong to $\hat{Q},$ where 0 is the zero element of $l_1$ (otherwise, $g^*(y_0)$ is not the minimum value in \eqref{D1}).  Due to Hahn-Banach separation theorem (see, e.g., \cite{DS}, Section V.2) there exists a sequence $(\k^j,\lambda^j)\in\reals\times l^\infty$ (where $\lambda^j=(\lambda_1^j,\lambda_2^j,\dots)$) such that
		\begin{equation}\label{E1}
		\begin{aligned}
			&\k^j\left(g^*(y_0)-\frac{1}{j}\right)+\delta^j \leq \inf_{(\theta,x)\in\hat{Q}}\{ \k^j\theta+\sum_{i=1}^{\infty}\lambda_i^jx_i \} \\
			&=\inf_{\gamma\in\mathcal{P}(G)}\{\k^j\theta+\int_{G}[\alpha(\psi_{\lambda^j}(f(y,u))-\psi_{\lambda^j}(y)) + (1-\alpha)(\psi_{\lambda^j}
		 (y_0)-\psi_{\lambda^j}(y))]\gamma(dy,du) \\
				&s.t. \, \theta\geq\int_{G}g(y,u)\gamma(dy,du)\},
			\end{aligned}
			\end{equation}
		where $\delta^j>0$ for all $j$ and $\psi_{\lambda^j}:=\sum_{i=1}^{\infty}\lambda_i^j\phi_i$. From the last formula it is easy to see that $\k^j\ge 0$. Let us show that, actually, $\k^j>0$. Indeed, if it was not the case and $\k^j=0$, then we would have
		\begin{equation*}
		\begin{aligned}
0< \delta^j\le &\min_{\gamma\in\mathcal{P}(G)}\int_{G}[\alpha(\psi_{\lambda^j}(f(y,u))-\psi_{\lambda^j}(y)) + (1-\alpha)(\psi_{\lambda^j}
		 (y_0)-\psi_{\lambda^j}(y))]\gamma(dy,du)\\
		=&	\min_{(y,u)\in G}\{(\alpha(\psi_{\lambda^j}(f(y,u))-\psi_{\lambda^j}(y)) + (1-\alpha)(\psi_{\lambda^j}
		 (y_0)-\psi_{\lambda^j}(y))\},
			\end{aligned}
			\end{equation*}
	which is a contradiction to \eqref{H3}. Thus, $\k^j>0$. Dividing \eqref{E1}  by $\k_j$, we obtain
		\begin{equation*}
		\begin{aligned}
		g^*(y_0)-\frac{1}{j}&	<\min_{\gamma\in \P(G)}\{\int_{G}\big(g(y,u)+{1\o \k^j}(\alpha(\psi_{\lambda^j}(f(y,u))-\psi_{\lambda^j}(y)) \\
		&+ (1-\alpha)(\psi_{\lambda^j}
		 (y_0)-\psi_{\lambda^j}(y))\big)\gm\}\\
		=\min_{(y,u)\in G}&\{g(y,u)+{1\o \k^j}(\alpha(\psi_{\lambda^j}(f(y,u))-\psi_{\lambda^j}(y)) + (1-\alpha)(\psi_{\lambda^j}
		 (y_0)-\psi_{\lambda^j}(y)))\}\\
		\le \mu^*_C(y_0)&.
			\end{aligned}
			\end{equation*}
Therefore, $g^*(y_0)\le\mu^*_C(y_0)$. Taking into account that, due to Proposition \ref{P5}, $g^*(y_0)=\mu^*(y_0)$, we obtain $\mu^*(y_0)\leq\mu^*_C(y_0)$. This, along with (\ref{e-inequality}), proves the (\ref{e:muC=mu}). \hfill{$\Box$}

\bigskip

A function $\psi\in D_N$ will be called a {\it solution of the $N$-approximating problem (\ref{eq-maxmin-N})} if
\begin{equation}\label{eq-maxmin-N-sol}
\begin{aligned}
\min_{(y,u)\in G}\{g(y,u)+\a(\psi(f(y,u))-\psi(y))+(1-\a)(\psi(y_0)-\psi(y))\}=\mu_N^*(y_0).
\end{aligned}
\end{equation}
From Proposition \ref{P6} it follows that, if $\psi\in D_N$ is a solution of (\ref{eq-maxmin-N}), then it solves the max-min problem (\ref{B5}) approximately in the sense that (compare with (\ref{B5-2})
$$
\min_{y\in Y}\{H_{\psi}(y)-(1-\a)(\psi(y)-\psi(y_0))\}=\mu^*(y_0)- \kappa_N,
$$
where $\kappa_N :=\mu^*(y_0)- \mu_N^*(y_0)\rightarrow 0$ as $N\rightarrow\infty$. Below we establish that a solution of the
$N$-approximating problem exists for any $N$ under a readily verifiable controllability-type assumption.

Let $\R_{y_0}$ be the reachable set for system \eqref{A1} in finite time. That is,
$$
\R_{y_0}:=\{y\ |\ y=y(t) \ \ {\rm  for\ some\ solution\ of\ (\ref{A1})}\ {\rm and\ for\ some}\ t\in \T  \}.
$$

\begin{Proposition}\label{P3.2}
Assume that
\begin{equation}\label{D7}
{\rm int}\,({\rm cl}\,\R_{y_0})\neq \emptyset.
\end{equation}
Then a solution of the $N$-approximating problem exists for any $N$.
\end{Proposition}

{\bf Proof.} Let us show first that the only function $\psi\in D_N$ satisfying the inequality
\begin{equation}\label{B11}
\alpha ({\psi}(f(y,u)) - {\psi}(y)) + (1-\alpha) ({\psi}(y_0) - {\psi}(y)) \geq 0 \,\, \ \ \ \forall \ (y,u)\in G,
\end{equation}
is $\psi\equiv 0$. Indeed, rewrite this inequality as
$$
\alpha {\psi}(f(y,u)) -{\psi}(y)+ (1-\alpha) {\psi}(y_0) \geq 0 \,\,  \ \ \ \forall \ (y,u)\in G.
$$
Let $(y(\cdot),u(\cdot))$ be an admissible process in \eqref{A1}. It follows from this inequality and \eqref{E4} that each term in \eqref{E4} is equal to zero, that is,
\begin{equation*}
\begin{aligned}
0&=\a\psi(f(y(t),u(t))-\psi(y(t))+(1-\a)\psi(y_0)\\
&=\a(\psi(y(t+1))-\psi(y(t)))+(1-\a)(\psi(y_0)-\psi(y(t))).
\end{aligned}
\end{equation*}
Hence,
$$
\psi(y(t+1))-\psi(y(t))={1-\a\o \a}(\psi(y(t))-\psi(y_0)) \quad\ \forall t\in\T.
$$
As can be readily verified, the latter implies $\psi(y(t))\equiv \psi(y_0)$ on any admissible trajectory. That is, $\psi(y)= \psi(y_0)$ for any $y\in \R_{y_0}$ and, consequently, $\psi(y)= \psi(y_0) \ \forall \ y\in {\rm cl}\,\R_{y_0}$. From  \eqref{D7} it now follows  that $\psi\equiv 0$ (due to
 \eqref{E2}).

Next, consider a maximizing sequence in the $N$-approximating problem, that is, for $k=1,2,\dots$ let $v^k=(v_1^k,\dots,v_N^k)\in\reals^N$ be such that the function
	\begin{equation*}
	\psi^k(y) := \sum_{i=1}^N v_i^k \phi_i(y)
	\end{equation*}
	satisfies the inequality
	\begin{equation}
	\label{eqn:psi k satisfies}
	g(y,u) + \alpha (\psi^k(f(y,u)) - \psi^k(y)) + (1-\alpha)(\psi^k(y_0) - \psi^k(y))\geq \mu^*_N(y_0)-\frac{1}{k}\ \ \ \forall \ (y,u)\in G.
	\end{equation}
		We will show that the sequence $v^k$, $k=1,2,\dots$ is bounded and, therefore, has a convergent subsequence. Assume to the contrary that
	there exists a subsequence $v^{k'}$ such that
	\begin{equation*}
	\label{eqn:vk unbounded}
	\lim_{k'\rightarrow\infty} |v^{k'}| = \infty,\,\, \lim_{k'\rightarrow\infty}\frac{v^{k'}}{|v^{k'}|} := \tilde{v},\,\, |\tilde{v}| = 1.
	\end{equation*}
	
	Dividing  (\ref{eqn:psi k satisfies}) by $|v^k|$ and passing to the limit along the subsequence $\{k'\}$, we obtain
	\begin{equation*}
	\alpha (\tilde{\psi}(f(y,u)) - \tilde{\psi}(y)) + (1-\alpha) (\tilde{\psi}(y_0) - \tilde{\psi}(y)) \geq 0 \,\,  \ \ \ \forall \ (y,u)\in G,
	\end{equation*}
	where
	\begin{equation*}
	\tilde{\psi}(y) := \sum_{i=1}^N \tilde{v}_i \phi_i(y).
	\end{equation*}
	We proved above that we must have $\tilde \psi\equiv 0$, which contradicts linear independence of the functions $\{\phi_i\}_{i=1}^N$.
	Thus, the sequence $v^k$ is bounded and there exists a subsequence such that $\disp 	\lim_{k'\rightarrow\infty} v^{k'} := v^*$.
	Passing to the limit in (\ref{eqn:psi k satisfies}) along this subsequence, we obtain
	\begin{equation*}
	g(y,u) + \alpha (\psi^*(f(y,u)) - \psi^*(y)) + (1-\alpha) (\psi^*(y_0) - \psi^*(y))\geq \mu^*_N(y_0)\,\, \ \ \ \forall \ (y,u)\in G,
	\end{equation*}
	where
	\begin{equation*}
	\psi^*(y) := \sum_{i=1}^N v_i^* \phi_i(y).
	\end{equation*}
	Therefore, $\psi^*(y)$ is a solution of the $N$-approximating problem. \hfill{$\Box$}

\section{Construction of near optimal controls}\label{Section-construction}


Let  \eqref{D7} hold true and let $\psi^N$ stand for a solution of the $N$-approximating problem. Motivated by \eqref{A16}, define the control $u^N$ by the equation
\begin{equation}\label{F8}
u^N(y)={\rm argmin}_{u\in A(y)}\{g(y,u)+\a\psi^N(f(y,u))\}
\end{equation}
and denote by $y^N$ the solution of the system
\begin{equation}\label{G9}
y^N(t+1)=f(y^N(t),u^N(y^N(t)))
\end{equation}
that satisfies the initial condition $y^N(0)=y_0$.

The next theorem (which is the main result of this section) states that, under appropriate conditions, $u^N(\cdot)$ and $y^N(\cdot)$ converge to the optimal control and the optimal trajectory as $N\to \infty$.

\begin{Theorem}\label{T1}
In addition to \eqref{D7} assume that the functions $f$ and $g$ are Lipschitz continuous and that the optimal solution $\g^*$ of problem \eqref{D1} is unique. Assume also that there exists an optimal admissible process $(\by(\cdot),\bu(\cdot)) $ such that:


{\rm (a)} For any $t\in \T$ there exists an open ball $Q_t$ centered at $\by(t)$ such that the minimizer $u^N(y)$ in the right hand side of \eqref{F8} is uniquely defined for $y\in Q_t$;

{\rm (b)} $u^N(\cdot)$ is Lipschitz continuous on $Q_t$ with Lipschitz constant independent of $N$ and $t$;

{\rm (c)} $y^N(t)\in Q_t\ \forall t\in \T$ for sufficiently large $N$.


Then
\begin{equation}\label{F7}
\begin{aligned}
&\lim_{N\to \infty} u^N(y^N(t))=\bu(t) \ \ \ \forall\ t\in\T,\\
&\lim_{N\to \infty}y^N(t)=\by(t)\ \ \ \forall\ t\in\T,\\
&\lim_{N\to \infty} V^N(y_0)=V(y_0),
\end{aligned}
\end{equation}
where $\disp V^N(y_0)=\sum_{t=0}^{\infty} \a^t g(y^N(t),u^N(y^N(t)))$.
\end{Theorem}

\bigskip

The proof of Theorem \ref{T1} is given at the end of this section. It is based on several propositions stated and proved below.

Consider the semi-infinite dimensional LP problem
\begin{equation}\label{F1}
\min_{\g\in W_N(y_0)} \int_{G} g(y,u)\g(dy,du)=:g^*_N(y_0),
\end{equation}
where
\begin{equation}\label{e-WN}
\begin{aligned}
&W_N(y_0):=\{\g\in \P(G)|\\
&\int_{G}\big[ \a(\phi_i(f(y,u))-\phi_i(y))
+(1-\a)(\phi_i(y_0)-\phi_i(y))\big]\gm=0,\quad i=1,\dots,N\},
\end{aligned}
\end{equation}
with $\phi_i$ as in Section \ref{Section-approx-max-min}.
Note that the set $W_N(y_0)$ is not empty (since $W_N(y_0)\supset W(y_0)\ \forall N$) and that it is compact in weak$^*$ topology. Therefore, the minimum in problem \eqref{F1} is reached. Note also that problem (\ref{F1}) is related to the $N$-approximating problem (\ref{eq-maxmin-N}). The latter is, in fact, dual  to the former, and the duality relationships  include, in particular,  the equality of the optimal values of these two problems, as established by the following lemma.


\begin{Lemma}\label{Lemma-duality-N}
The optimal value of  problem (\ref{F1}) is equal to the optimal value of the  $N$-approximating problem (\ref{eq-maxmin-N}):
\begin{equation}\label{e-duality-N}
\disp g_N^*(y_0)=\mu_N^*(y_0).
\end{equation}
\end{Lemma}

{\bf Proof.}
 Denote
\begin{equation*}\label{A7}
\mu_N(\psi,y_0):=\min_{(y,u)\in G}\{g(y,u)+\a (\psi(f(y,u))-\psi(y))+(1-\a)(\psi(y_0)-\psi(y))\}\ \ \ \forall \ \psi\in D_N,
\end{equation*}
so that $\disp \mu_N^*(y_0)=\sup_{\psi\in D_N}\mu(\psi,y_0)$. For any $\g\in W_N(y_0)$ and $\psi\in D_N$, we have
\begin{equation*}
\begin{aligned}
\mu_N(\psi,y_0)\le& \int_{G}(g(y,u)+\a (\psi(f(y,u))-\psi(y))+(1-\a)(\psi(y_0)-\psi(y)))\gm\\
=&\int_{G}g(y,u)\gm,
\end{aligned}
\end{equation*}
which implies that
\begin{equation}\label{K6}
\mu_N^*(y_0)\le g_N^*(y_0).
\end{equation}
The proof of the opposite inequality is similar to the proof of Lemma \ref{e-duality-inf}. Namely,
we first show  that, for any
$\psi\in D_N$, the inequality $\disp\min_{(y,u)\in G}\{\a(\psi(f(y,u))-\psi(y))+(1-\a)(\psi(y_0)-\psi(y))\}\le 0$ is valid (compare with \eqref{H3}). Then we introduce the set $\hat{Q}\subset\reals\times \reals^N$,
		\begin{equation*}
				\begin{aligned}
			\hat{Q} := &\{(\theta,x) |\, \theta\geq\int_{G}g(y,u)\gamma(dy,du),\, x=(x_1,\dots,x_N),  \\
			&x_i=\int_{G}[\alpha(\phi_i\left(f\left(y,u\right)\right)-\phi_i\left(y\right)) + \left(1-\alpha\right)(\phi_i\left(y_0\right)-\phi_i\left(y\right))]\gamma\left(dy,du\right),\\
			&i=1,\dots,N,\,\gamma\in\mathcal{P}(G)\}
		\end{aligned}
		\end{equation*}
and use convex separation theorem to separate $\hat Q$ from the point $(g^*_N(y_0)-1/j,0)$ ($0$ being the zero element of $\reals^N $; see (\ref{E1})). This will lead to the inequality $g^*_N(y_0)\le \mu^*_N(y_0)$, which, along with \eqref{K6}, will prove the validity of (\ref{e-duality-N}).
\hfill{$\Box$}

\begin{Proposition}\label{P4.1}
The following relations are true:

{\rm(i)} $\disp \lim_{N\to \infty} \rho_H(W_N(y_0),W(y_0))=0$;

{\rm(ii)} $\disp \lim_{N\to \infty} g_N^*(y_0)=(1-\a)V(y_0)$;

{\rm(iii)} If the optimal solution $\g^*$ of problem \eqref{D1} is unique, then $\disp \lim_{N\to \infty}\g^N=\g^*$, where $\g^N$ is any optimal solution of \eqref{F1}.


\end{Proposition}

{\bf Proof.} Since $W(y_0)\subset W_N(y_0)$, to prove (i), it suffices to show that
$$\disp \lim_{N\to \infty}\sup_{\g\in W_N}\rho(\g,W)=0.$$
 Assume this is not true. Then there exist a sequence $\{\g_{N'}\}\in W_{N'}$ and a number $\d>0$ such that $\rho(\g_{N'},W)\ge \d$ for all $N'$. From weak$^*$ compactness of $\P(G)$ it follows that there exists $\bg\in \P(G)$ and a subsequence of $\{\g_{N'}\}$ (we do not relabel) such that
$$
\lim_{N'\to \infty}\g^{N'}=\bg,
$$
therefore,
\begin{equation}\label{G1}
\rho(\bg,W)\ge \d.
\end{equation}
On the other hand, since $\g_{N'}\in W_{N'}$, we have
$$
\int_{G}\big[\a(\phi_i(f(y,u))-\phi_i(y))+(1-\a)(\phi_i(y_0)-\phi_i(y))\big]\g_{N'}(dy,du)=0
$$
for any $i$ (provided that it less or equal than $N'$). Hence,
$$
\int_{G}\big[ \a(\phi_i(f(y,u))-\phi_i(y))+(1-\a)(\phi_i(y_0)-\phi_i(y))\big]\bg(dy,du)=0 \ \ \ \forall \ i.
$$
Due to the approximating property of $\{\phi_i \} $, it implies that
$$
\int_{G}\big[ \a(\phi(f(y,u))-\phi(y))+(1-\a)(\phi(y_0)-\phi(y))\big]\bg(dy,du)=0 \ \ \ \forall \phi\in C(Y)
$$
and, thus, it follows that
   $\bg\in W$. This contradicts \eqref{G1} and completes the proof of statement (i). Note that (i) implies that
$\disp \lim_{N\to \infty} g_N^*(y_0)=g^*(y_0)$, which, due to \eqref{D3}, proves the validity of statement (ii).
The validity of statement (iii) follows from the fact that, due to (i), any partial limit of an optimal solution of \eqref{F1} is an optimal solution of \eqref{D1}.
\hfill{$\Box$}

\bigskip

\begin{Proposition}\label{P4.3-0}
Among the optimal solutions of problem \eqref{F1}, there exists one (denoted below as $\g^N$) that is presented as a convex combination of at most $N+1$  Dirac measures with concentration points in $G$. More precisely,
\begin{equation}\label{F2}
\g^N=\sum_{j=1}^{K_N}\b_j^N\d_{(y_j^N,u_j^N)}, \ \ \ {\rm where} \ \ \ \b_j^N>0, \ j=1,\dots,K_N\le N+1, \ \ \ \ \disp \sum_{j=1}^{K_N}\b_j^N=1
\end{equation}
and where \ $\d_{(y_j^N,u_j^N)}$ are the Dirac measures concentrated at $(y_j^N,u_j^N)\in G$.
 Moreover, the concentration points $(y_j^N,u_j^N), \ j=1,\dots,K_N\ $ satisfy the following relationships:
\begin{equation}\label{G5}
\begin{aligned}
&u_j^N={\rm argmin}_{u\in A(y_j^N)} \{g(y^N_j,u)+\a\psi^N(f(y^N_j,u))\},\\
&y_j^N={\rm argmin}_{y\in Y}\{H_{\psi^N}(y)-(1-\a)\psi^N(y)\},
\end{aligned}
\end{equation}
where $\psi^N$ is a solution of the $N$-approximating problem (\ref{eq-maxmin-N}).
\end{Proposition}

{\bf Proof.}
Let $W_N^*$ be the optimal solution set in \eqref{F1} and let $\g^N$ be one of its extreme points. (The set of extreme points of $W_N^*$ is not empty; see, e.g., \cite{DS}, Section V.8). Being an extreme point of $W_N^*$,  $\g^N$ must be an extreme point of the set $W_N$ as well.  Any extreme point of  $W_N$ can be represented as a sum of  at most $N+1$ Dirac measures (see, e.g., Theorem A.5 in \cite{Rubio}). Thus, the representation \eqref{F2} is valid.

Let us prove the validity of \eqref{G5}. Taking into account the fact that $\g^N$ is an optimal solution of \eqref{F1} and that $\g^N\in W_N(y_0)$ we have
\begin{equation*}\label{G2}
\begin{aligned}
g_N^*&(y_0)=\int_G g(y,u)\g^N(dy,du)\\
&=\int_G [g(y,u)+\a(\psi^N(f(y,u))-\psi^N(y))+(1-\a)(\psi^N(y_0)-\psi^N(y))]\g^N(dy,du).
\end{aligned}
\end{equation*}
 Using \eqref{F2}, we obtain
\begin{equation}\label{G3}
\begin{aligned}
g_N^*(y_0)= \sum_{j=1}^{K_N}\b_j^N[g(y_j^N,u_j^N)+&\a(\psi^N(f(y_j^N,u_j^N))-\psi^N(y_j^N))\\
&+(1-\a)(\psi^N(y_0)-\psi^N(y_j^N))].
\end{aligned}
\end{equation}
On the other hand, due to Lemma \ref{Lemma-duality-N} and due to the fact that $\psi^N$ is a solution of the $N$-approximating problem (see (\ref{eq-maxmin-N-sol})), we have
\begin{equation}\label{G4}
\begin{aligned}
g_N^*(y_0)=\mu_N^*(y_0)=\min_{(y,u)\in G}\{g(y,u)+&\a(\psi^N(f(y,u))-\psi^N(y))\\
&+(1-\a)(\psi^N(y_0)-\psi^N(y))\}.
\end{aligned}
\end{equation}
Comparing \eqref{G3} and \eqref{G4} and taking into account that  $\b_j^N>0$, we obtain that, for all $j=1,
\dots,K_N$,
\begin{equation*}
\begin{aligned}
&g(y_j^N,u_j^N)+\a(\psi^N(f(y_j^N,u_j^N))-\psi^N(y_j^N))+(1-\a)(\psi^N(y_0)-\psi^N(y_j^N))=\\
&\min_{(y,u)\in G}\{g(y,u)+\a(\psi^N(f(y,u))-\psi^N(y))+(1-\a)(\psi^N(y_0)-\psi^N(y))\},
\end{aligned}
\end{equation*}
which is equivalent to \eqref{G5}. \hfill{$\Box$}

\begin{Proposition}\label{P4.3}
Let $(\by(\cdot),\bu(\cdot))$ be an optimal process in \eqref{A1} such that the conditions (a),(b) and (c) of Theorem \ref{T1} are satisfied and let $\g^N$ be an optimal solution of \eqref{F1} that is  represented in the form \eqref{F2}. Then, for any $t$, there exist points
$(u^N_{j_N},y^N_{j_N})\in\{(y^N_{j},u^N_{j}),\,j=1,\dots,K_N\}$ such that
\begin{equation}\label{F5}
(\by(t),\bu(t))=\lim_{N\to \infty}(y^N_{j_N},u^N_{j_N}).
\end{equation}
\end{Proposition}
{\bf Proof.}
Since the optimal solution $\g^*$  of the IDLP problem \eqref{D1} is unique and since the discounted
occupational measure $\g_{\bar u} $ generated by the optimal control $\bu(\cdot)$ must be an optimal solution of   \eqref{D1}
(due to Proposition \ref{P5}; see also (\ref{A3}) and (\ref{D1})), one may conclude that $\g_{\bar u}=\g^*$.
That is,  $\g^*$ is the discounted occupational measure generated by $\bu(\cdot)$. Due to the definition of the latter (see \eqref{EE6}) one comes to the conclusion that
\begin{equation}\label{F4}
\g^*(B_r(\by(t),\bu(t)))>0\ \ \ \ \forall \ \ t\in\T \ \ \ {\rm and} \ \ \ \ \forall \ r>0,
\end{equation}
where $B_r(\by,\bu):=\{(y,u)|\,|y-\by|+|u-\bu|<r\}$.

Assume that the statement of the proposition is false. Then there exist $r>0$, $\ t\in\T$ and a sequence $N_i\to \infty$ such that
$$
(y_j^{N_i},u_j^{N_i})\notin B_r (\by(t),\bu(t))\ \ \ \forall \ j=1,\dots,K_{N_i}.
$$
By \eqref{F2}, it follows that
\begin{equation*}
\g^{N_i}(B_r(\by(t),\bu(t)))=0.
\end{equation*}
Due to statement (iii) of Proposition \ref{P4.1},  $\disp \lim_{N_i\to \infty}\g^{N_i}=\g^*$. This implies
$$
0=\lim_{i\to \infty}\g^{N_i}(B_r(\by(t),\bu(t)))\ge\g^*(B_r(\by(t),\bu(t)))
$$
(due to the semicontinuity property of probability measures; see Theorem 2.1 in \cite{Bill}).
The latter contradicts \eqref{F4}. The obtained contradiction proves the proposition.
\hfill{$\Box$}

\bigskip

{\bf Proof of Theorem \ref{T1}.} Let $t\in \T$ and let  $(u^N_{j_N},y^N_{j_N})$ be as in \eqref{F5}. Comparing formulas \eqref{F8} and \eqref{G5}, we conclude that
$$
u^N_{j_N}=u^N(y^N_{j_N}).
$$
Therefore,
\begin{equation}\label{F6}
|\bu(t)-u^N(\by(t))|\le |\bu(t)-u^N_{j_N}|+|u^N(y^N_{j_N})-u^N(\by(t))|.
\end{equation}
Due to \eqref{F5} and due to the assumed continuity of the map $u^N(\cdot)$ in a neighbourhood of $\by(t) $ (see assumptions (b) and (c) of Theorem \ref{T1}), we have
\begin{equation}\label{F9}
\lim_{N\to \infty} u^N(\by(t))=\bu(t) \quad\hbox{for all }t.
\end{equation}

Subtracting the equation $\disp \by(t+1)=f(\by(t),\bu(t))$ from the equation $\disp y^N(t+1)=f(y^N(t),u^N(y^N(t)))$  and taking into account the fact that the
functions $f(\cdot,\cdot)$ and $u^N(\cdot)$ are  Lipschitz continuous, we obtain
\begin{equation*}
\begin{aligned}
|y^N(t+1)-&\by(t+1)|\le L_1(|y^N(t)-\by(t)|+|u^N(y^N(t))-\bu(t)|)\\&\le L_1(|y^N(t)-\by(t)|+|u^N(y^N(t))-u^N(\by(t))|+|u^N(\by(t))-\bu(t)|)\\
&\le L_2(|y^N(t)-\by(t)|+|u^N(\by(t))-\bu(t)|)
\end{aligned}
\end{equation*}
were  $L_1$ and $L_2$ are some appropriately chosen positive constants.
The inequalities above (along with \eqref{F9} and the discrete time analog of Gronwall-Bellman lemma) allow to conclude that
$$
|y^N(t)-\by(t)|\to 0\quad \hbox{as }N\to \infty.
$$
Also,
$$
|u^N(y^N(t))-\bu(t)|\le |u^N(y^N(t))-u^N(\by(t))|+|u^N(\by(t))-\bu(t)|\to 0 \hbox{ as }N\to \infty.
$$
Thus, the first two relationships in \eqref{F7} are proved. The validity of  the third relationship in \eqref{F7} follows from the first two and from the fact that $g $ is Lpischitz continuous. \hfill{$\Box$}

\section{Numerical example}\label{Section-example}

An optimal solution (\ref{F2}) of the semi-infinite LP problem (\ref{F1}) as well as  its optimal value $g^*_N $ and an optimal solution of the $N$-approximating problem (\ref{eq-maxmin-N})   can be found numerically using, e.g., the algorithm discussed in \cite{GRT}. Also, once an optimal solution of the $N$-approximating problem  is found, one can  construct a control $u^N(y)$ as a minimizer in (\ref{F8}), the latter  being near optimal in (\ref{A111}) for $N$ large enough (as has been established by Theorem \ref{T1}).

 Note that conditions under which the statements of Theorems \ref{T1}  are valid are difficult to verify. Their verification, however, is not needed for the construction of the control $u^N(y) $. After this control is constructed, one can find the admissible process $(y^N(\cdot), u^N(\cdot)) $ (from \eqref{G9}) and subsequently find the corresponding value of the objective function $V_N(y_0)$. Since
$$
g_N^*(y_0)=\mu_N^*(y_0)\le \mu^*(y_0)=(1-\a)V(y_0)\le (1-\a)V_N(y_0),
$$
the difference  $V_N(y_0)-V(y_0)$ is less than or equal to the difference $V_N(y_0)- g_N(y_0)(1-\a)^{-1} $, the latter provides us with a \lq\lq measure of  near optimality" of the found control.

Let us illustrate the way a near optimal control can be constructed  with the help of  a numerical example.

{\bf Example 1.} Consider the optimal control problem
(\ref{A111}) with
$$
y=(y_1,y_2),\ \ \  u=(u_1,u_2), \ \ \ g(y,u) =  -y_1(t)u_2(t) + y_2(t)u_1(t)
$$
and with $\ f(y,u)=(f_1(y,u), f_2(y,u)) $, where
\begin{equation}
\label{ex1:f}
f_i(y,u) = \frac{1}{2}y_i - \frac{1}{2}u_i, \ \ \ \ \ i=1,2.
\end{equation}
Let the map $U(y)$ be constant-valued:
$\ U(y)=U:=[-1,1]\times [-1,1] $ and let $\ Y:=[-1,1]\times [-1,1]$. One can readily verify that $\ f(y,u)\in Y\ \ \forall u\in U $ if $\ \ y\in Y $ and, therefore, $\ A(y)=U$
(see (\ref{e-A})).
The semi-infinite LP problem (\ref{F1}) was formulated for this problem with the use of the monomials $\phi_{i_1,i_2}(y)=y_1^{i_1}y_2^{i_2}, \ i_1,i_2= 0,1,..., \mathcal{I}, $ as the functions $\phi_i(\cdot) $ defining $W_N(y_0)$ in (\ref{e-WN}). Note that the number of constraints $N$ in (\ref{e-WN}) is equal to $(\mathcal{I}+1)^2$ in this case.
This problem and the corresponding $N$-approximating problem were solved numerically   with the use of the algorithm similar to one described in \cite{GRT} for $\mathcal{I}=7 $ ($N=64$). The discount factor was taken to be equal to $0.9$ ($\alpha = 0.9$), and the initial conditions were taken to be as follows:
\begin{equation}
\label{ex1:ic}
y_1(0)= \frac{1}{2}, \ \ \ \ \ \ y_2(0) = \frac{1}{4}.
\end{equation}
In particular, the coefficients  of the expansion $\lambda_{i_1,i_2}^N $ defining the optimal solution of the  $N$-approximating problem,
$$
\psi^N(y_1,y_2)= \sum_{0\leq i_1+i_2\leq \mathcal{I}} \lambda_{i_1,i_2}^N y_1^{i_1}y_2^{i_2}
$$
were found, and the optimal value of the semi-infinite LP problem was evaluated to be approximately equal to $-1.013 $ ($g_N^*(y_0)\approx- 1.013$).


\medskip

For $(y_1, y_2) =(y_1(0), y_2(0))= (0.5, 0.25) $ (see (\ref{ex1:ic})), the minimizer of the problem
\begin{equation}\label{e-nonlinear-optim}
\min_{(u_1,u_2)\in U}\{ -y_1u_2 + y_2u_1 + 0.9\psi^N(\frac{1}{2}y_1-\frac{1}{2}u_1, \frac{1}{2}y_2-\frac{1}{2}u_2 )\}
\end{equation}
was numerically found (using MATLAB) to be equal to $(-1.000, 1.000) $. Thus, we take $u^N(0)= (-1,1) $, which after substitution into the equations of the dynamics, gives
$$
y^N_1(1)= \frac{1}{2}(0.5) - \frac{1}{2}(-1)= 0.75, \ \ \ \ y^N_2(1)= \frac{1}{2}(0.25) - \frac{1}{2}(1)= -0.375
$$
For $(y_1, y_2) =(y_1^N(1), y_2^N(1))$, the minimizer of the problem (\ref{e-nonlinear-optim}) was found to be $\approx (-0.552, 1.000) $ and we take $u^N(1) =(-0.552, 1.000) $. The latter being substituted into the equations of the dynamics allows one to obtain
$$
y^N_1(2)= \frac{1}{2}(0.75) - \frac{1}{2}(-0.552)= 0.651, \ \ \ \ y^N_2(2)= \frac{1}{2}(-0.375) - \frac{1}{2}(1)\approx -0.688.
$$
This process has been repeated $50$ times,  and the results of the first $10$ time steps are shown in the table below.

\begin{center}
\begin{tabular}{|c|c|c|c|c|}
	\hline
	$t$ & $y_1^N(t)$ & $y_2^N(t)$ & $u_1^N(t)$ & $u_2^N(t)$ \\
	\hline
	0 &0.500 &0.250 & -1.000 & 1.000\\
	\hline
	1& 0.750 & -0.375 & -0.552 & 1.000 \\
	\hline
	2 & 0.651 & -0.688 & 1.000 & 1.000 \\
	\hline
	3 & -0.174 &-0.844 & 1.000 & 1.000 \\
	\hline
	4 & -0.587 & -0.922 & 1.000 & -1.000 \\
	\hline
5 & -0.794 & 0.039 & 1.000 & -1.000 \\
\hline
6 & -0.897 & 0.520 & -1.000 & -1.000 \\
\hline
7 & 0.052 & 0.760 & -1.000 & -1.000 \\
\hline
8& 0.526 & 0.880 & -1.000 & 1.000 \\
\hline
9& 0.763 & -0.060 & -1.000 & 1.000 \\
\hline
10& 0.881 & -0.530 & 1.000 & 1.000 \\

	 \hline
\end{tabular}
\end{center}

\begin{figure}
	\hspace{0.25\textwidth}
	\includegraphics[width=0.5\textwidth]{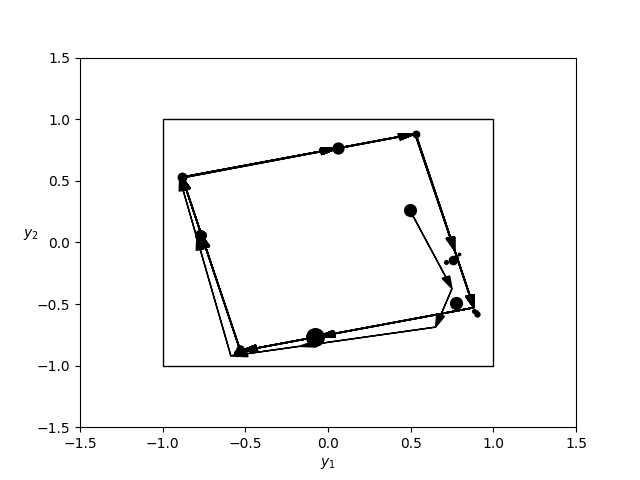}
	\caption{The state trajectory -  50 time steps.}
	\label{fig:ex1_trajectory_duals}
\end{figure}

Note that, starting from the moment $t=2$, the controls take only values of $1$ or $-1$. Also, starting from this moment, the sequence of controls appeared to be periodic with the period $T=8 $ (that is, $\ (u^N_1(t+8), u^N_2(t+8))= (u^N_1(t), u^N_2(t)) \ \forall \ t\geq 2 $).
The obtained state trajectory appears to be converging to a \lq\lq square like figure" as shown in
 Fig \ref{fig:ex1_trajectory_duals}. The concentration points of the Dirac measures in the expansion  (\ref{F2}) (these and the corresponding weights were found as a part of an optimal solution of the semi-infinite LP problem; see Section 4 in \cite{GRT}) are marked  with dots in Fig \ref{fig:ex1_trajectory_duals}, the size of which are scaled proportionally to the magnitude of their respective weights. Having in mind Proposition \ref{P4.3}, one can expect that the  optimal state trajectory must come close to at least some of these dots, and, as one can see,
  the obtained state trajectory has a similar property  (passing near or just going through  these dots).

  The value of the objective function  thus obtained was evaluated to be approximately equal to $-9.972$  (${V}^{N}(y_0)\approx -9.972$). Consequently,
  $$|{V}^{N}(y_0)-g_N^*(y_0)(1-\a)^{-1}|\approx - 9.972 - (-10.13) = 0.158,$$
   which indicates that the value of the objective function obtained with the use of the constructed control is within a close proximity of the optimal one.




{\bf Remark.} Note that the control $u^N(y)$  defined as a minimizer in (\ref{F8}) is near optimal only in a neighbourhood of the optimal trajectory that starts at $y_0$. If, instead of $\psi^N(\cdot) $,
an  approximation of the optimal value function is used in the right-hand-side of (\ref{F8}), then the corresponding minimizer will be near optimal for all $y\in Y$. A most common approach to finding an approximation of the optimal value function is based on a discretization of the state space and solving the corresponding dynamic programming equation for thus obtained finite state space process (see, e.g., Appendix A in \cite{BCD}). The finer is  the grid of the  discretization, the better are approximations of the optimal value and  the optimal control. In contrast to this approach, the function $\psi^N(\cdot) $, being an approximate solution of the max-min problem (\ref{B5}), is obtained via solving the semi-infinite dimensional (SID) linear programming problem (\ref{F1}), and the proximity of an approximation depends on the number $N$ of the test functions in (\ref{e-WN}). While solving the SIDLP problem (\ref{F1}) requires a discretization of the state space, the grid of the discretization does not need to be fine. In fact, the algorithm proposed in \cite{GRT} aims at finding the grid points (their number being no more than $N$) that are concentrated around the optimal trajectory. Arguably, this may require less computational efforts than the classic approach although, of course, more research is needed to validate this claim.

\section{Generalization of Theorem \ref{T1} and a heuristic numerical algorithm}\label{heuristic}

Finding  a minimizer in (\ref{F8}) may  be a challenging task since the problem in the right-hand-side of (\ref{F8}) is generally  not of the convex programming class (even in the case when $g(y,u)$ is convex and $f(y,u)$ is linear in $u$). This difficulty can be avoided since the only property of the minimizer $u^N(y)$ that was used in proving that it is near optimal for large $N$ is that it satisfies the equalities
\begin{equation}\label{e-essential-equalities}
u^N(y^N_j) = u^N_j, \ \ j=1,..., K_N,
\end{equation}
where $(y^N_j, u^N_j)$ are the concentration points of the Dirac measures in the presentation (\ref{F2}) (see Proposition \ref{P4.3-0}). Having this in mind, one can establish  a  more general result suggesting that a simpler way of constructing a near optimal control is possible. To state the result, let us introduce the following assumption.

ASSUMPTION II:
Let the set of the concentration points
$\{ (y^N_j, u^N_j), \ j=1,...,K_N\} $  does not contain points with the same $y$-coordinates and different $u$-coordinates.

\medskip

Note that, under Assumption II (which is assumed to be valid everywhere in this section), one can define the function $\xi^N(\cdot)$,
\begin{equation}\label{e-finite-function}
\xi^N(y^N_j): = u^N_j, \ \ j=1,..., K_N.
\end{equation}
This function  is defined on the set of $y$-coordinates of the concentration points $\{ (y^N_j, u^N_j), \ j=1,...,K_N\} $ and it takes values in the set of  $u$-coordinates of these points.

\medskip

We will call  a function $u^N(\cdot)$  an extension of $\xi^N(\cdot)$ onto $Y$ if $u^N(y)\in A(y) \ \ \forall \ y\in Y$ and if the equalities (\ref{e-essential-equalities}) are satisfied. (As has been mentioned above, the   minimizer in (\ref{F8}) is an example of such  an extension.) The following theorem establishes that, under conditions similar to those of Theorem \ref{T1}, an extension of $\xi^N(\cdot)$ onto $Y$ is near optimal for $N$ is large enough.

\begin{Theorem}\label{T1-1} Let $u^N(\cdot)$ be an extension of $\xi^N(\cdot)$ onto $Y$ and let $y^N(\cdot)$  be the solution of (\ref{G9}) satisfying the initial condition $y^N(0)=y_0$. Let \eqref{D7} be satisfied,
the functions $f$ and $g$ be Lipschitz continuous and  the optimal solution $\g^*$ of problem \eqref{D1} be unique. Assume also that there exists an optimal admissible process $(\by(\cdot),\bu(\cdot)) $ such that:


{\rm (a)} For any $t\in \T$ there exists an open ball $Q_t$ centered at $\by(t)$ such that
 $u^N(\cdot)$ is Lipschitz continuous on $Q_t$ with Lipschitz constant independent of $N$ and $t$;

{\rm (b)} $y^N(t)\in Q_t\ \forall t\in \T$ for sufficiently large $N$.


Then the relationships (\ref{F7}) are valid.

\end{Theorem}

{\bf Proof.} The proof follows exactly the same steps as those in the proof of Theorem \ref{T1}.
\hfill{$\Box$}



  In the case when the map $A(y)$
  is constant-valued, it is natural do define an extension $u^N(\cdot) $ of $\xi^N(\cdot)$ in such way that it takes constant values on  some neighbourhoods of the $y$-components of the concentration points
$\ \{ (y^N_j, u^N_j), \ j=1,...,K_N\} $. That is, $\ u^N(y)= u^N_j $ for all $y$ in some \lq\lq sufficiently small" neighborhood of $\ y^N_j$ for all  $j=1,...,K_N$.
 (Note that, due to Proposition \ref{P4.3}, one may expect that the optimal admissible trajectory $\by(\cdot)$ is contained in the union of these neighbourhoods   if $N$ is large enough.)

Using the idea of such piecewise constant extension of $\xi^N(\cdot)$,
 one can propose the following heuristic algorithm for the construction of a control that can be a \lq\lq good candidate" for being near optimal in (\ref{A111}):

 \medskip

 {\bf Heuristic algorithm for the construction of a near optimal control.}

 1. Find an optimal solution  of the semi-infinite LP problem (\ref{F1}). That is, find the concentration points $\ \{ (y^N_j, u^N_j)\}$ of the Dirac measures and the weights $\{\beta^N_j\}$, $ \ j=1,...,K_N,$ in the presentation (\ref{F2}), and evaluate
 the optimal value $g_N^*(y_0) $ of this problem.

 \medskip

 2.  Among the $y$-components of the concentration points, choose one which is  closest to $y_0$. That is, choose
 $\ j_0:= {\rm argmin}_j\{|y^N_j  - y_0|\}$ and define
  \begin{equation*}\label{e-step-1}
  u^N(0):=u^N_{j_0}, \ \ \ \ \   y^N(1):= f(y_0, u^N(0)).
 \end{equation*}
 Assume that $ u^N(0), ..., u^N(t-1) $ and $ y^N(1), ...,  y^N(t) $, where
 $t\geq 1  $, have been defined. Then $ u^N(t) $ and $  y^N(t+1)$ are defined as follows:
 \begin{equation*}\label{e-step-t}
  u^N(t):=u^N_{j_t}, \ \ \ \ \   y^N(t+1):= f( y^N(t),  u^N(t)),
 \end{equation*}
 where $j_t:={\rm argmin}_j\{|y^N_j -  y^N(t)|\} $.

 \medskip

3. Continue this process until the moment $t=T $ when the $\ \frac{\a^{T+1}}{1-\a}$ becomes small enough, making the finite sum $\sum_{t=0}^{T} \a^t g( y^N(t),  u^N(t)):={V}^{N}(y_0)$  a good approximation for the value of the objective function on the infinite time horizon. Evaluate the difference $\ |{V}^{N}(y_0)-g_N^*(y_0)(1-\alpha)^{-1}|$ that provides a measure of near optimality of the constructed control.

{\bf Remark.}
{\rm If for some $t\geq 0 $  the minimizer in the problem $\ {\rm min}_j\{|y^N_j -  y^N(t)|\} $ is not unique, then one can take $\ j_t:= {\rm argmax}_{j\in J_t }\{\beta^N_j\} $, where
  $\ J_t:={\rm Argmin}_j\{|y^N_j -  y^N(t)|\}$. That is, among the concentration points the $y$-components of which are equally close to $ y^N(t)$, one can choose one that has the  greatest weight. Note that a
  weight
   $\beta^N_j$ in the sum $\sum_{j=1}^{K_N}\b_j^N\d_{(y_j^N,u_j^N)}$ (see (\ref{F2})) can be interpreted as an approximation of the number of times  the optimal trajectory  \lq\lq attends" a vicinity of the corresponding $y^N_j $ (each subsequent attendance being time discounted).

   }

Let us illustrate the way the algorithm outlined above works with the help of  the  example considered in Section \ref{Section-example}.

{\bf Example 1 (continued).} A sample of the concentration points $\ y_j^N= (y_{1,j}^N, y_{2,j}^N)$, $\ u_j^N= (u_{1,j}^N, u_{2,j}^N)$ of the Dirac measures entering (\ref{F2}) and the values of the corresponding weights $\{\beta^N_j\}$, $ \ j=1,...,K_N$   is shown in the table below (these being a part of the solution of  the SIDLP problem of Example 1).

\begin{center}
\begin{tabular}{|c|c|c|c|c|}
	\hline
	$y_{1,j}^N$ & $y_{2,j}^N$ & $u_{1,j}^N$ & $u_{2,j}^N$ & $\beta_j^N$ \\
	\hline
	-0.5375 & -0.875 & 1 & -1 & 0.0913 \\
	\hline
	0.5 & 0.25 & -1 & 1 & 0.0820 \\
	\hline
	0.7625 & -0.4875 & 1 & 1 & 0.0642 \\
	\hline
	-0.775 & 0.05 & 1 & -1 & 0.0538 \\
	\hline
	0.8875 & -0.5625 & 1 & 1 & 0.0536 \\
	\hline
-0.0875 & -0.75 & 1 & 1 & 0.0525 \\
	 \hline
\end{tabular}
\end{center}
(The concentration points  with weights that are less than  $10^{-2} $ were discarded, the number of terms in (\ref{F2}) after the discarding being equal to $27$.)

\begin{figure}
	\hspace{0.25\textwidth}
	\includegraphics[width=0.5\textwidth]{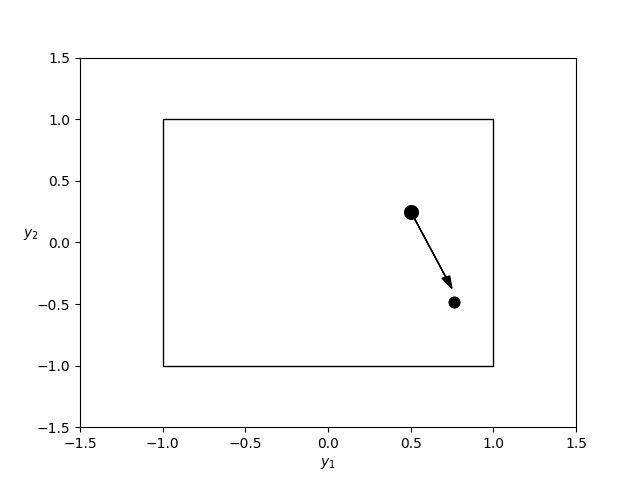}
	\caption{\ The  state trajectory - time step 1.}
	\label{fig:ex1_1-iter}
\end{figure}

\begin{figure}
	\hspace{0.25\textwidth}
	\includegraphics[width=0.5\textwidth]{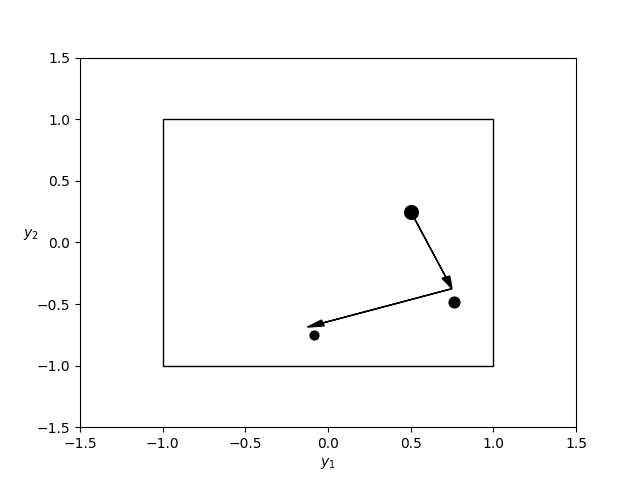}
	\caption{The state trajectory - time steps 1 and 2.}
	\label{fig:ex1_2-iter}
\end{figure}

Among the concentration points shown in the table above, there is one with the $y$-coordinates $(0.5, 0.25)$,  these coinciding with the initial conditions (\ref{ex1:ic}). Thus, we take ${u}^N(0) $ to be equal to the corresponding $u $-coordinates of this point, ${u}^N(0)=(-1,1)  $, which being substituted into the equations of the dynamics, leads to ${y}^N(1)= (0.75, -0.375)$. The closest to ${y}^N(1)$ (of all $y$ components of the concentration points) is the pair $(0.7625, -0.4875)$; the latter is marked by a dot in Fig. \ref{fig:ex1_1-iter}. The corresponding $u$-components of the concentration point are $(1,1)$. We, therefore, take  ${u}^N(1)=(1,1)  $ and  obtain ${y}^N(2)= (-0.125, -0.6875)$. The closest to ${y}^N(2)$ are the $y$-coordinates $(-0.0875, -0.75) $ marked in Fig \ref{fig:ex1_2-iter}. The corresponding $u$-components of the concentration point are $(1,1)$, and, consequently,  ${u}^N(2)=(1,1) $, etc.

\begin{figure}
	\hspace{0.25\textwidth}
	\includegraphics[width=0.5\textwidth]{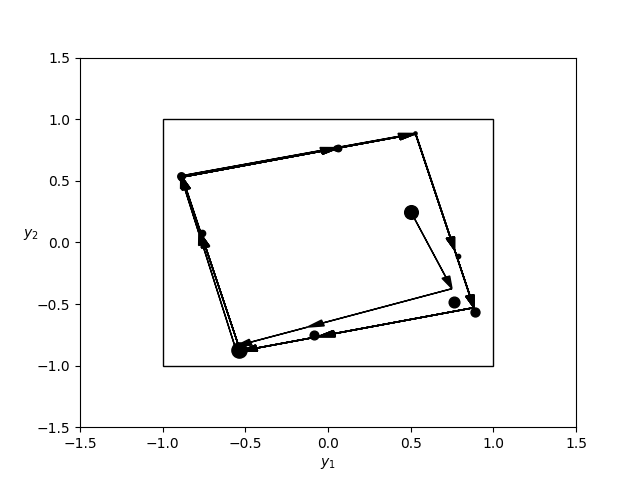}
	\caption{The state trajectory -  50 time steps.}
	\label{fig:ex1_trajectory}
\end{figure}

This process has been repeated $50$ times and the resulted state trajectory is depicted in Fig. \ref{fig:ex1_trajectory}. Note that it looks similar to that depicted in Fig \ref{fig:ex1_trajectory_duals}.
 The value of the objective function  thus obtained was evaluated to be approximately equal to $-10.03$  (${V}^{N}(y_0)\approx-1.003$). Finally, we obtain
 $$|{V}^{N}(y_0)-g_N^*(y_0)(1-\a)^{-1}|\approx (-10.03) - (-10.13)= 0.1,$$
  the error being less than in the case when the control defined as a minimizer   (\ref{F8}) was used.




{\bf Remark.} As has been mentioned above, results obtained in this paper are in many ways analogous to similar results in continuous time setting obtained in \cite{GRT}. 
Let us briefly outline  relationships between these two sets of results. Firstly, all results of the present paper are established for the case when the control set may depend on the state variables, with the only restriction on this dependence  being that the map $U(\cdot) $ is upper semicontinuous and compact valued. The consideration in \cite{GRT} was for the fixed control set, that is without the dependence of the latter on the state variables. Moreover, allowing such  dependence in the continuous time setting would make the consideration  much more technical (requiring, e.g., additional regularity assumptions about the map $U(\cdot)$). Secondly, the $max$ in the max-min problem (\ref{B5})  is over the class of bounded lower semicontinuous functions, this being in contrast to the continuous time max-min counterpart of (\ref{B5}), in which $max$ (or rather $sup$) was over the set of smooth functions. While  maximizer in (\ref{B5}) always exists, a maximizer in the continuous counterpart of (\ref{B5}) may not exist and its existence needs to be assumed, this being a restrictive assumption
(compare  Theorem \ref{P2.1} and Proposition 2.1 in \cite{GRT}).

As far as the algorithmic part of the paper is concerned, finding  a minimizer in problem (\ref{F8}) is generally a more difficult task than finding a minimizer in the continuous time counterpart of (\ref{F8}) (see Section 3 in \cite{GRT}). For example, the latter is of convex programming while the former is generally  not in case  $g(y,u)$ is convex and $f(y,u)$ is linear in $u$. In Section 6, we proposed a way allowing one to avoid solving (\ref{F8}), this being also possible in continuous time setting (such an opportunity was not investigated in \cite{GRT}).

\section*{Acknowledgments} A significant part of the work on the paper was conducted during I. Shvartsman's  extended visit to the Department of Mathematics at Macquarie University (Sydney, Australia). The research was supported by the ARC Discovery Grants DP130104432 and DP150100618 and by a Macquarie University Start-Up Grant.


\begin{thebibliography}{99}

\bibitem{Adelman} (MR2125136) [10.1287/moor.1040.0109]
\newblock D. Adelman and D. Klabjan,
\newblock \doititle{Duality and existence of optimal policies in generalized joint replenishment},
\newblock \emph{Mathematics of Operations Research}, {\textbf 30} (2005), 28--50.

\bibitem{Ash} (MR0435321)
 \newblock  R. Ash,
 \newblock \emph{Measure, Integration and Functional Analysis},
 \newblock Academic Press, 1972.

\bibitem{Aub} (MR1134779)
\newblock J.-P. Aubin,
\newblock \emph{Viability Theory},
\newblock Birkh\"auser, 1991.

\bibitem{BCD} (MR1484411) [10.1007/978-0-8176-4755-1]
\newblock M. Bardi and I. Capuzzo-Dolcetta,
\newblock \emph{Optimal Control and Viscosity Solutions of Hamilton-Jacobi-Bellman Equations, Systems and Control: Foundations and Applications,}
\newblock Birkh\"auser, Boston, 1997.

\bibitem{Bert}  (MR3644954)
\newblock D. Bertsekas, 
\newblock \emph{Dynamic Programming and Optimal Control,}
\newblock Athena Scientific, Belmont, MA, 2017.


\bibitem{BhBo} (MR1411505) [10.1214/aop/1065725192]
\newblock A. G. Bhatt and V. S. Borkar,
\newblock \doititle{Occupation measures for controlled Markov processes: Characterization and optimality},
\newblock \emph{Annals of Probability}, {\textbf {24}} (1996), 1531--1562.

\bibitem{Bill}
P. Billingsley,
\emph{Convergence of Probability Measures}, John Wiley \& Sons, New York, 1968.



\bibitem{B12} (MR2918256)
\newblock J. Blot,
\newblock A Pontryagin principle for infinite-horizon problems under constraints,
\newblock \emph{Dynamics of Continuous, Discrete and Impulsive Systems Series B: Applications and Algorithms}, {\textbf {19}} (2012), 267--275.

\bibitem{Vivek} (MR0950347) [10.1007/BF00353877]
\newblock V. S. Borkar,
\newblock \doititle{A convex analytic approach to Markov decision processes},
\newblock \emph{Probability Theory and Related Fields,} {\textbf {78}} (1988), 583--602.

\bibitem{BGQ} (MR2772196) [10.1007/s00245-010-9120-y]
\newblock R. Buckdahn, D. Goreac and M. Quincampoix,
\newblock \doititle{Stochastic optimal control and linear programming approach},
\newblock \emph{Appl. Math. Optim.} {\textbf {63}} (2011), 257--276.

\bibitem{CHL} (MR3155340) [10.1007/978-3-642-76755-5]
\newblock D. A. Carlson, A. B. Haurier and A. Leizarowicz,
\newblock \emph{Infinite Horizon Optimal Control. Deterministic and Stochastic Processes},
\newblock Springer, Berlin, 1991.

\bibitem{DS} (MR1009162)
\newblock N. Dunford and J. T. Schwartz,
\newblock \emph{Linear Operators, Part I, General Theory},
\newblock Wiley \& Sons, Inc., New York, 1988.

\bibitem{Gai-F-Leb} (MR2421325) [10.1137/060676398]
\newblock L. Finlay, V. Gaitsgory and I. Lebedev,
\newblock \doititle{Duality in linear programming problems related to deterministic long run average problems of optimal control,}
\newblock \emph{SIAM J. Control and Optimization}, {\textbf {47}} (2008), 1667--1700.

\bibitem{F-V} (MR1009341) [10.1137/0327060]
\newblock W. H. Fleming and D. Vermes,
\newblock \doititle{Convex duality approach to the optimal control of diffusions,}
\newblock \emph{SIAM J. Control Optimization}, {\textbf {27}} (1989), 1136--1155.

\bibitem{Gai8} (MR2082704) [10.1137/S0363012903424186]
\newblock V. Gaitsgory,
\newblock \doititle{On representation of the limit occupational measures set of control systems with applications to singularly perturbed control systems},
\newblock \emph{SIAM J. Control and Optimization}, {\textbf {43}} (2004), 325--340.

\bibitem{GPS-0} [0.1109/CDC.2016.7798950]
\newblock V. Gaitsgory, A. Parkinson and I. Shvartsman,
\newblock \doititle{Linear programming formulation of a discrete time infinite horizon optimal control problem with time discounting criterion},
\newblock \emph{Proceedings of 55th IEEE  Conference on Decision and Control (CDC)}, 2016,  Las Vegas, USA.


\bibitem{GPS-1} (MR3693843)
\newblock V. Gaitsgory, A. Parkinson and I. Shvartsman,
\newblock \doititle{Linear Programming Formulations of Deterministic Infinite Horizon Optimal Control Problems in Discrete Time}, \newblock \emph{Discrete and Continuous Dynamical Systems, Series B},  {\textbf {22(10)} (2017), 3821--3838.


\bibitem{GQ} (MR2556353) [10.1137/070696209]
\newblock V. Gaitsgory and M. Quincampoix,
\newblock \doititle{Linear programming approach to deterministic infinite horizon optimal control problems with discounting},
\newblock \emph{SIAM J. Control and Optim.}, {\textbf {48}} (2009), 2480--2512.

\bibitem{GQ-1} (MR3057046) [10.1016/j.na.2013.03.015]
\newblock V. Gaitsgory and M. Quincampoix,
\newblock \doititle{On sets of occupational measures generated by a deterministic control system on an infinite time horizon,}
\newblock \emph{Nonlinear Analysis (Theory, Methods} \& \emph{Applications)}, {\textbf {88}} (2013), 27--41.

\bibitem{GR} (MR2248173) [10.1137/040616802]
\newblock V. Gaitsgory and S. Rossomakhine,
\newblock \doititle{Linear programming approach to deterministic long run average problems of optimal control,}
\newblock \emph{SIAM J. of Control and Optimization}, {\textbf {44}} (2006), 2006--2037.

\bibitem{GRT} (MR2918249)
\newblock V. Gaitsgory, S. Rossomakhine and N. Thatcher,
\newblock \doititle{Approximate solution of the HJB inequality related to the infinite horizon optimal control problem with discounting},
\newblock \emph{Dynamics of Continuous, Discrete and Impulsive Systems
Series B: Applications and Algorithms}, {\textbf{19}} (2012), 65--92.

\bibitem{Goreac-Serea} (MR3041666) [10.1051/cocv/2011183]
\newblock D. Goreac and O.-S. Serea,
\newblock \doititle{Linearization techniques for $L^{\infty} $ - control problems and dynamic programming principles in classical and $L^{\infty} $ control problems},
\newblock \emph{ESAIM: Control, Optimization and Calculus of Variations}, {\textbf {18}} (2012), 836--855.

\bibitem{GruneSIAM98} (MR1626864) [10.1137/S0363012997315919]
\newblock L. Gr\"une,
\newblock \doititle{Asymptotic controllability and exponential stabilization of nonlinear control systems at singular points},
\newblock \emph{SIAM J. Control Optim.}, {\textbf {36}} (1998), 1495--1503.

\bibitem{GruneJDE98} (MR1637521) [10.1006/jdeq.1998.3451]
\newblock L. Gr\"une,
\newblock \doititle{On the relation between discounted and average optimal value functions},
\newblock \emph{J. Diff. Equations}, {\textbf {148}} (1998), 65--69.

\bibitem{Her-Her-Lasserre} (MR1404981)
\newblock D. Hernandez-Hernandez, O. Hernandez-Lerma and M. Taksar,
\newblock {The linear programming approach to deterministic optimal control problems},
\newblock \emph{Appl. Math.}, {\textbf {24}} (1996), 17--33.

\bibitem{Jean} O. Hernandez-Lerma and J.B. Lasserre, \lq\lq The linear Programmimg Approach", in the volume {\it Handbook of Markov Decision Processes: Methods and Applications},  Edited by E.A. Feinberg and A. Shwartz, Springer, 2012.

\bibitem{Kamo}
\newblock A. Kamoutsi, T. Sutter, P. M. Esfahani, and J. Lygeros,
\newblock {On Infinite Linear Programming and the Moment Approach to Deterministic Infinite Horizon
Discounted Optimal Control Problems},
\newblock \emph{IEEE Control System Letters}, {\textbf {1}} (2017), 134--139.

\bibitem{Adelman-1} (MR2348232) [10.1287/moor.1070.0252]
\newblock D. Klabjan and D. Adelman,
\newblock \doititle{An Infinite-dimensional linear programming algorithm for deterministic semi-Markov decision processes on Borel spaces},
\newblock \emph{Mathematics of Operations Research}, {\textbf {32}} (2007), 528--550.

\bibitem{Kurtz} (MR1616514) [10.1137/S0363012995295516]
\newblock T. G. Kurtz and R. H. Stockbridge,
\newblock \doititle{Existence of Markov controls and characterization of optimal Markov controls},
\newblock \emph{SIAM J. on Control and Optimization}, {\textbf {36}} (1998), 609--653.

\bibitem{Lass-Trelat} (MR2421324) [10.1137/070685051]
\newblock J. B. Lasserre, D. Henrion, C. Prieur and E. Tr\'elat,
\newblock \doititle{Nonlinear optimal control via occupation measures and LMI-relaxations},
\newblock \emph{SIAM J. Control Optim.}, {\textbf {47}} (2008), 1643--1666.


\bibitem{QS} (MR2580138) [10.1016/j.na.2009.11.024]
\newblock M. Quincampoix and O. Serea,
\newblock \doititle{The problem of optimal control with reflection studied through a linear optimization problem stated on occupational measures},
\newblock \emph{Nonlinear Anal.} {\textbf {72}} (2010), 2803--2815.


\bibitem{Rubio} (MR0832189)
\newblock J. E. Rubio,
\newblock \emph{Control and Optimization. The Linear Treatment of Nonlinear Problems},
\newblock Manchester University Press, Manchester, 1986.

\bibitem{Stockbridge} (MR1043943) [10.1214/aop/1176990944]
\newblock R. H. Stockbridge,
\newblock \doititle{Time-Average control of a martingale problem. Existence of a stationary solution},
\newblock \emph{Annals of Probability}, {\textbf {18}} (1990), 190--205.

\bibitem{Stockbridge1} (MR1043944) [10.1214/aop/1176990945]
\newblock R. H. Stockbridge,
\newblock \doititle{Time-average control of a martingale problem: A linear programming formulation},
\newblock \emph{Annals of Probability}, {\textbf {18}} (1990), 206--217.

\bibitem{Filar92} (MR1189274) [10.1007/BF00939913]
\newblock R. Sznajder and J. A. Filar,
\newblock \doititle{Some comments on a theorem of Hardy and Littlewood},
\newblock \emph{J. Optimization Theory and Applications}, {\textbf {75}} (1992), 201--208.

\bibitem{Vinter} (MR1205987) [10.1137/0331024]
\newblock R. Vinter,
\newblock \doititle{Convex duality and nonlinear optimal control},
\newblock \emph{SIAM J. Control and Optim.,} {\textbf {31}} (1993), 518--538.

\bibitem{Z14} (MR3309874) [10.1007/978-3-319-08034-5]
\newblock A. Zaslavski,
\newblock \emph{Stability of the Turnpike Phenomenon in Discrete-Time Optimal Control Problems},
\newblock Springer, 2014.

\bibitem{Z06} (MR2164615)
\newblock A. Zaslavski,
\newblock \emph{Turnpike Properties in the Calculus of Variations and Optimal Control,}
\newblock Springer, New York, 2006.

\bibitem{Z06a} (MR3306947) [10.1007/978-3-319-08828-0]
\newblock A. Zaslavski,
\newblock \emph{Turnpike Phenomenon and Infinite Horizon Optimal Control,}
\newblock Springer Optimization and Its Applications, New York, 2014.



}


\end{thebibliography}
\end{document}